%% file: main.tex
\documentclass[nopreprintline]{elsarticle}
\usepackage[margin=1in]{geometry}
\usepackage{lineno}
\input{config}

\usepackage{etoolbox}
\newcommand*\linenomathpatch[1]{%
  \cspreto{#1}{\linenomath}%
  \cspreto{#1*}{\linenomath}%
  \csappto{end#1}{\endlinenomath}%
  \csappto{end#1*}{\endlinenomath}%
}
\linenomathpatch{abstract}
\linenomathpatch{equation}
\linenomathpatch{gather}
\linenomathpatch{align}

\def\appendixname{\!\!}
\newcommand{\apdxref}[1]{\hyperref[#1]{Appendix~\ref{#1}}}

\author[a1]{Joseph~Hart}
\ead{joshart@sandia.gov}
\affiliation[a1]{
   organization={Scientific~Machine~Learning,~Sandia~National~Laboratories},
   addressline={ 1611~Innovation~Pkwy.~SE},
   city={Albuquerque,~NM},
   postcode={87123},
   country={United~States~of~America},
}
   \affiliation[a2]{
   organization={Department~of~Mathematics,~Brigham~Young~University},
   addressline={275~TMCB},
   city={Provo,~UT},
   postcode={84602},
   country={United~States~of~America},
}
\cortext[cor1]{Corresponding author (\href{mailto:joshart@sandia.gov}{joshart@sandia.gov}).}

\author[a1,a2]{Shane~A.~McQuarrie}
\ead{smcquar@sandia.gov}

\author[a1]{Zachary~Morrow}
\ead{zbmorro@sandia.gov}

\author[a1]{Bart~van~Bloemen~Waanders}
\ead{bartv@sandia.gov}

\journal{Journal of Computational Physics}

\date{}

\begin{document}


\thispagestyle{empty}

\begin{abstract}
\input{abstract}
\end{abstract}


\maketitle

\input{body}

\appendix
\bigskip\noindent\textbf{Appendices}
\input{appendices}

\bibliographystyle{abbrv}
\bibliography{references}

\end{document}

%% file: config.tex
\usepackage{amsfonts,amssymb,amsmath,amsthm,amsbsy}
\usepackage{comment}
\usepackage{hyperref}
\usepackage{graphicx}
\usepackage{empheq}
\usepackage{subfigure}
\usepackage{microtype}
\usepackage[dvipsnames]{xcolor}
\usepackage[normalem]{ulem}  
\usepackage{url}
\usepackage[nameinlink]{cleveref}

\definecolor{blue}{RGB}{0, 61, 165}
\definecolor{royalblue}{RGB}{35, 87, 178}
\definecolor{sandiablue}{RGB}{0, 50, 90}
\definecolor{sandiared}{RGB}{130, 36, 51}
\definecolor{sandiapurple}{RGB}{148,0,211}
\definecolor{burntorange}{cmyk}{0, 0.65, 1, .09}

\hypersetup{
    colorlinks=true,
    linkcolor=blue,
    citecolor=blue,
    urlcolor=blue,
    pdftitle={Toward real-time optimization through model reduction and model discrepancy sensitivities},
    pdfauthor={Joseph Hart, Shane A. McQuarrie, Zachary Morrow, Bart van Bloemen Waanders}
}

\crefformat{equation}{#2(#1)#3}
\numberwithin{equation}{section}

\theoremstyle{definition}

\theoremstyle{plain}

\theoremstyle{definition}

\renewcommand{\vec}[1] {\ensuremath{\boldsymbol{#1}}}
\def\u{{\vec{u}}}
\def\x{{\vec{x}}}
\def\z{{\vec{z}}}
\def\y{{\vec{y}}}
\def\q{{\vec{q}}}

\def\f{{\vec{f}}}

\def\s{{\vec{s}}}
\def\R{{\mathbb{R}}}
\def\V{{\vec{V}}}
\def\v{{\vec{v}}}
\def\A{{\vec{A}}}
\def\M{{\vec{M}}}

\def\B{{\vec{B}}}

\def\S{{\vec{S}}}

\def\xivec{{\vec{\xi}}}

\def\H{{\vec{H}}}
\def\P{{\vec{P}}}
\def\U{{\vec{U}}}
\def\W{{\vec{W}}}
\def\d{{\vec{\delta}}}
\def\t{{\vec{\theta}}}
\def\ztilde{{\tilde{\z}}}
\newcommand{\trp}{^{\mathsf{T}}}

\usepackage[section]{algorithm}
\usepackage[noend]{algpseudocode}
\algnewcommand{\LineComment}[1]{\State \textcolor{gray}{\texttt{\# #1}}}
\algrenewcomment[1]{\hfill\textcolor{gray}{\texttt{\# #1}}}

\title{\Large\textbf{Toward real-time optimization through model reduction \\and model discrepancy sensitivities}}

\usepackage{fancyhdr}   

%% file: abstract.tex
Optimization problems arise in a range of scenarios, from optimal control to model parameter estimation. In many applications, such as the development of digital twins, it is essential to solve these optimization problems within wall-clock-time limitations. However, this is often unattainable for complex systems, such as those modeled by nonlinear partial differential equations. One strategy for mitigating this issue is to construct a reduced-order model (ROM) that enables more rapid optimization. In particular, the use of nonintrusive ROMs---those that do not require access to the full-order model at evaluation time---is popular because they facilitate optimization solutions can be computed within the wall-clock-time requirements. However, the optimization solution will be unreliable if the iterates move outside the ROM training data. This article proposes the use of hyper-differential sensitivity analysis with respect to model discrepancy (HDSA-MD) as a computationally efficient tool to augment ROM-constrained optimization and improve its reliability. The proposed approach consists of two phases: (i) an offline phase where several full-order model evaluations are computed to train the ROM, and (ii) an online phase where a ROM-constrained optimization problem is solved, $N=\mathcal{O}(1)$ full-order model evaluations are computed, and HDSA-MD is used to enhance the optimization solution using the full-order model data. Numerical results are demonstrated for two examples, atmospheric contaminant control and wildfire ignition location estimation, in which a ROM is trained offline using inaccurate atmospheric data. The HDSA-MD update yields a significant improvement in the ROM-constrained optimization solution using only one full-order model evaluation online with corrected atmospheric data.

%% file: body.tex
\section{Introduction}
\label{sec:intro}

Advancements in modeling and computing capabilities have opened up the possibility of creating digital twins for complex systems. This
can impact a range of applications
including additive manufacturing, fusion energy, aerospace systems, and wildfire containment. 
Stated briefly, a digital twin is an goal-oriented 
coupling of physical and virtual assets with a bidirectional information
flow between the assets~\cite{NAP26894}.
Central to the digital twin concept is the ability to adapt and inform decision-making as data is collected
 in real-time---which is defined relative to the timescale of the data
 collection and decisions. This requires solving optimization problems
 constrained by numerical models to (i) calibrate the digital twin to physical system
 data, and (ii) use the digital twin to inform decisions, such as actions to control the physical system.
However, the computational cost of model
 evaluations makes real-time optimization unattainable in many
 applications. This article seeks to enable real-time optimization
 for systems modeled by partial differential equations (PDEs)
  through a novel combination of
 non-intrusive model reduction and hyper-differential sensitivity
 analysis with respect to model discrepancy.

 Model reduction seeks to alleviate the computational burden of solving a
 high-dimensional system,  referred to as the full-order
 model (FOM), by identifying and exploiting structure to
 build a low-dimensional system, referred to as the reduced-order
 model (ROM). Classical model reduction is intrusive in the sense that
 it learns low-dimensional structure from training data but also requires
 access to the FOM to construct and/or solve the ROM. This is unsuitable for many applications requiring real-time
 computation because the cost of FOM access inside the
 optimization loop remains prohibitively expensive. Non-intrusive model
 reduction methods seek to overcome this limitation by learning a
 low-dimensional system exclusively from training data so that
 additional access to the FOM is not needed. This ensures
 that the computational cost of reduced-order model evaluations is
 independent of the FOM dimension. However, this gain in
 computational efficiency generally results in a loss in accuracy, especially when
 predicting outside of the regime represented by the training data. 
 This presents a challenge for optimization problems constrained by the ROM, as the optimizer's iterations may drift away from the training data regime, leading to an unreliable solution.

 To improve the reliability of ROM-constrained optimization, we propose
 to augment it with hyper-differential sensitivity analysis with respect to model discrepancy, herein referred
 to as HDSA-MD. The approach originated
 in~\cite{hart2023hdsaoptimal,hart2024hdsasampling} to facilitate
 optimization for problems where a high-fidelity model can only be
 evaluated a small number of times. Specifically, these articles
 consider optimization constrained by a low-fidelity model, derived
 from physics simplifications such as omission of small-magnitude
 forces or linearization, and use $N=\mathcal O(1)$ evaluations of a high-fidelity
 PDE model to update the low-fidelity
 optimization solution. The benefits of HDSA-MD are that high-fidelity
 derivative information is not required and the evaluations can be
 performed outside of the optimization loop. This article combines
 HDSA-MD and non-intrusive reduced-order modeling to enable real-time
 optimization. Computational efficiency is achieved using
 ROM-constrained optimization that is independent of the FOM, and reliability is improved via the HDSA-MD optimal solution update.

We demonstrate our approach by training an operator inference ROM to
enable real-time atmospheric contaminant control and by training a
neural operator ROM to enable real-time wildfire ignition source
estimation. In both examples, we require that the ROMs are trained
``offline'' with inaccurate estimates of the atmospheric conditions.
These pre-trained ROMs are insufficient since the state of the
atmosphere at decision time is generally different than in the training data. However,
by augmenting the ROM-constrained optimization with HDSA-MD, we show
that using as little as one FOM evaluation with accurate
atmospheric data is sufficient to enhance the ROM-constrained optimization.

There is a vast literature on the use of ROMs to reduce the
computational cost of PDE-constrained optimization. Goal-oriented
model reduction in the service of optimization is introduced
in~\cite{buithanh2006goob}. To improve ROM
accuracy,~\cite{yue2013optlinearrom,qian2017rbopt,zahr2015pderom,hassdonk2023rbmlromopt}
propose various approaches that construct a sequence of optimization
problems and iteratively enrich the ROM. Specialized ROMs for Hessian
approximation are considered in~\cite{heinkenschloss2018rom}. ROMs are
used to enable shape optimization
in~\cite{antil2011DDandBT,antil2012ROM} and risk-averse optimization
in~\cite{heinkenschloss2018CVaR,zou2022CVaR}. Recent work has also explored
the use of derivative information to train neural operator ROMs in the service of solving optimization
problems~\cite{dipnet_2022,Luo_2025}. This article differs from previous works in
that we focus on real-time applications where a non-intrusive ROM is
trained offline in preparation for ROM-constrained optimization
online, yet a small number (as little as 1) of FOM
evaluations are feasible online to enhance the optimization
solution using HDSA-MD. 

The article is organized as follows. \Cref{sec:opt_form} presents the optimization problem under consideration and the corresponding ROM-constrained optimization problem derived from it. We present our HDSA-MD formulation for ROM-constrained optimization in \Cref{sec:hdsa}. \Cref{sec:atmosphere_transport} and~\Cref{sec:wildfire} demonstrate our approach for ROM-constrained optimization arising in control of atmospheric contaminants and estimation of wildfire ignition, respectively. Concluding remarks are made in \Cref{sec:conclusion}.

\section{Optimization formulation} \label{sec:opt_form}
Consider FOM-constrained optimization problems of the form
\begin{subequations}
\begin{empheq}[left={\texttt{FOMCO}\empheqbiglbrace\hspace{5mm}}]{align}
    \label{eqn:true_objective}
    & \displaystyle\min_{\z \in \R^{n_z}} J(\S(t;\z,\xivec),\z)
    \\ \nonumber
    & \textup{where } \S(t;\z,\xivec) \textup{ solves}
    \\ \label{eqn:true_constraint}
    &\begin{aligned}
    & \dot{\u}(t) = \f(\u(t),\z, \xivec) \qquad t \in [0,T]
    \\
    & \u(0) = \mathring{\u}.
    \end{aligned}
\end{empheq}
\end{subequations}

Time is denoted by $t \in [0,T]$, $T>0$ is the final time, $\u(t) \in
\R^{n_u}$ is the state, $\z \in \R^{n_z}$ corresponds to control, design, or
inversion parameters, $\xivec \in \R^{n_\xi}$ is a vector of model parameters,
$\f: \R^{n_u} \times \R^{n_z} \times \R^{n_\xi} \to \R^{n_u}$ defines
the system dynamics, $\mathring{\u} \in \R^{n_u}$ is the initial state, and
$\S(t;\z,\xivec)$ is the solution operator that maps $(\z,\xivec)$ to
the solution of~\eqref{eqn:true_constraint}. We assume that
$\S(t;\z,\xivec)$ is uniquely defined for each $(\z,\xivec)$. The
objective function $J:L^2([0,T]) \times \R^{n_z} \to \R$ depends on the solution operator $
\S$, and hence each objective function
evaluation requires solving~\eqref{eqn:true_constraint}. We seek to
solve~\texttt{FOMCO} with fixed model parameters $\xivec$; however,
the model parameters are uncertain, so variability of $\xivec$
plays an important role in our analysis. In our numerical results,
$\xivec$ corresponds to atmospheric winds, which are only known within
a short time window before the optimization solution is required.

The~\texttt{FOMCO} problem arises in many applications where a PDE is
discretized in space, yielding the state vector $\u(t) \in
\R^{n_u}$. In this case, $n_u$ is the number of nodes in the spatial discretization, making~\cref{eqn:true_constraint} a high-dimensional system.
Throughout the article we will refer to $\z \in \R^{n_z}$ as
the decision variable, noting that it may correspond to control, design, or
inversion parameters. We formulate~\texttt{FOMCO} with the state
written as a function of time, i.e., $\u(t)$, but the decision variables $\z$ and
model parameters $\xivec$ as vectors in Euclidean space. In general,
we may consider $\z$ and $\xivec$ that are functions of space and/or time, but
we assume that they are discretized in~\texttt{FOMCO} to simplify the
presentation. Where it is relevant, we will use weighted inner
products to bring information from the function spaces and
discretization into our analysis.

Assuming that~\texttt{FOMCO} cannot be solved within wall-clock-time
requirements, we consider a process in which FOM
evaluations are performed offline in preparation for online
computation to approximate the solution of~\texttt{FOMCO}.
In the offline phase, we generate data by evaluating the FOM for $M$
decision variables instances $\z_k \in \R^{n_z}$, $k=1,2,\dots,M$, to produce a set of state
trajectories $\u_k(t)$, $k=1,2,\dots,M$, satisfying~\eqref{eqn:true_constraint}. From these trajectories, a
low-dimensional linear subspace,
$\text{span} \{ \v_1,\v_2,\dots,\v_{n_{\hat{u}}} \}$, is determined,
where $\v_i \in \R^{n_u}$, $i=1,2,\dots,n_{\hat{u}}$, are basis
vectors computed from FOM data. We let
$\V = \begin{pmatrix} \v_1 & \v_2 & \cdots &
  \v_{n_{\hat{u}}} \end{pmatrix} \in \R^{n_u\times n_{\hat{u}}}$
denote the basis matrix and represent the state
$\u(t)$ by the time-dependent coordinates
$\widehat{\u}(t) =
(\hat{u}_{1}(t),\ldots,\hat{u}_{n_{\hat{u}}}(t))\in\mathbb{R}^{n_{\hat{u}}}$
as
\begin{align*}
    \u(t)
    \approx \V \widehat{\u}(t)
    = \sum_{i=1}^{n_{\hat{u}}}\hat{u}_{i}(t) \v_i.
\end{align*}
This reduces the state degrees-of-freedom from $n_u$ to $n_{\hat{u}}$,
where $n_{\hat{u}} \ll n_u$. The full dynamics $\f:\R^{n_u} \times \R^{n_z}
\times \R^{n_\xi} \to \R^{n_u}$ are approximated by reduced dynamics
$\widehat{\f}:\R^{n_{\hat{u}}} \times \R^{n_z} \times \R^{n_\xi} \to
\R^{n_{\hat{u}}}$.
There are a plurality of ways to construct $\widehat{\f}$. For instance, Galerkin projection defines $\widehat{\f}(\widehat{\u},\z,\xivec) = \V\trp \f(\V \widehat{\u},\z,\xivec)$. However, this is intrusive in the sense that evaluating $\widehat{\f}$ requires access to $\f$. This article focuses on non-intrusive methods---those which only require data output from FOM solves---to construct $\widehat{\f}$.  In particular, we demonstrate our framework using operator inference~\cite{peherstorfer2016opinf} and neural operator surrogates~\cite{kovachki2023neural}, which are detailed in~\Cref{sec:atmosphere_transport} and~\Cref{sec:wildfire}, respectively.

Due to the complexity of learning a non-intrusive ROM from limited FOM
data, we consider the case where the model parameters $\xivec$ are
fixed to a nominal value $\tilde{\xivec}$ for FOM data generation
and ROM training. We suppress the dependence of $\widehat{\f}$ on
$\xivec$ with the implication that $\widehat{\f}$ is trained using
data with $\xivec=\tilde{\xivec}$. We will account for uncertainty
in $\xivec$ through our proposed HDSA-MD approach.

 Given a basis matrix $\V$ and approximate dynamics $\widehat{\f}$, we consider the ROM-constrained optimization problem
 \begin{subequations}
 \begin{empheq}[left={\texttt{ROMCO}\empheqbiglbrace\hspace{5mm}}]{align}
    \label{eqn:rom_objective}
    & \displaystyle\min_{\z \in \R^{n_z}} J(\V \widehat{\S}(t;\z),\z)
    \\ \nonumber
    & \textup{where } \widehat{\S}(t;\z) \textup{ solves}
    \\ \label{eqn:rom_constraint}
    &\begin{aligned}
    & \dot{\widehat{\u}}(t) = \widehat{\f}(\widehat{\u}(t),\z)
    \qquad t \in [0,T]
    \\
    & \widehat{\u}(0) = \V\trp \mathring{\u}.
    \end{aligned}
\end{empheq}
\end{subequations}
Thanks to the dimension reduction $n_{\hat{u}} \ll n_u$, the ROM
solution operator $\widehat{\S}(t;\z)$ can be evaluated much faster
than the FOM solution operator $\S(t;\z,\xivec)$ and thereby enable
real-time performance.

Let $\ztilde$ denote the~\texttt{ROMCO} minimizer that is computed online, but that is insufficient due to the discrepancy between the FOM and ROM. This discrepancy stems from three sources: (i) restriction to the basis matrix $\V$; (ii) failing to learn the complete FOM dynamics; and (iii) misspecification of $\xivec$ in training. Our goal in this article is to improve the~\texttt{ROMCO} solution $\ztilde$ in real-time using $N$
FOM solves. In particular, we focus on the case where only $N =
\mathcal O(1)$ FOM evaluations are achievable within wall-clock-time
limitations. We assume that the ``true'' model parameters
$\xivec=\xivec^\dagger$ are available online so that the $N$ FOM
solves use the model parameters $\xivec^\dagger$ rather than the model
parameters $\tilde{\xivec}$ used to generate data for ROM
training. The goal is to update $\ztilde$ to better approximate
the~\texttt{FOMCO} solution when
$\xivec=\xivec^\dagger$.

\section{HDSA with respect to model discrepancy}
\label{sec:hdsa}

Consider the model discrepancy parameterized optimization problem
 \begin{subequations}
 \begin{empheq}[left={\texttt{ROMCO-MD}\empheqbiglbrace\hspace{5mm}}]{align}
    \label{eqn:rom_objective}
    & \displaystyle\min_{\z \in \R^{n_z}} J(\V \widehat{\S}(t;\z) + \d(t ;\z,\t),\z)
    \\ \nonumber
    & \textup{where } \widehat{\S}(t;\z) \textup{ solves}
    \\ \label{eqn:hdsa_constraint}
    &\begin{aligned}
    & \dot{\widehat{\u}}(t) = \widehat{\f}(\widehat{\u}(t),\z)
    \qquad t \in [0,T]
    \\
    & \widehat{\u}(0) = \V\trp \mathring{\u}.
    \end{aligned}
\end{empheq}
\end{subequations}
\texttt{ROMCO-MD} corresponds to endowing~\texttt{ROMCO} with a model discrepancy operator $\d(t;\z,\t)$, parameterized by $\t \in \R^{n_{\theta}}$, which seeks to account for the discrepancy between FOM and ROM, i.e., $\S(t;\z,\xivec^\dagger)-\V \widehat{\S}(t;\z)$. In what follows, we demonstrate how to use this model discrepancy endowment to enhance the~\texttt{ROMCO} solution.

\subsection{Model discrepancy calibration} \label{ssec:dis_calibration}
Assume that the FOM is evaluated $N$ times to provide
discrepancy data $\{ \z_\ell, \S(t;\z_\ell,\xivec^\dagger)-\V \widehat{\S}(t;\z_\ell)
\}_{\ell=1}^N$. We seek to parameterize $\d$ and calibrate its parameters $\t$ so that
\begin{align} \label{eqn:discrepancy_approx}
\d(t;\z_\ell,\t) \approx \S(t;\z_\ell,\xivec^\dagger)-\V \widehat{\S}(t;\z_\ell), \qquad \ell=1,2,\dots,N. \ 
\end{align}
In general, the discrepancy between FOM and ROM depends nonlinearly on $\z$. However, we cannot expect to learn this nonlinearity since the training dataset only contains $N=\mathcal O(1)$ FOM evaluations. Rather, the model discrepancy operator is parameterized as affine in $\z$:
 \begin{align}
\label{eqn:affine_dis}
\d(t; \z,\t) = \u_0(t;\t) + \vec{L}(t;\t) \M_\z \z.
\end{align}
Here, $\u_0:[0,T] \times \R^{n_\theta} \to \R^{n_u}$ and $\vec{L}:[0,T] \times \R^{n_\theta} \to \R^{n_u \times n_z}$ are the intercept and linear terms in the affine operator. The mass matrix $\M_\z \in \R^{n_z \times n_z}$, whose $(i,j)$ entry is the inner product of the $i^{th}$ and $j^{th}$ decision-variable basis function\footnote{The mass matrix $\M_\z$ is the identity matrix if $\z$ does not arise from discretization of a function.}, ensures that the discrepancy parameterization respects the structure of the decision variables if they arise from spatial and/or temporal discretization. The intercept $\u_0$ and linear operator $\vec{L}$ are parameterized as linear functions of $\t$ and hence $\d(t;\z,\t)$ depends linearly on $\t$.

To calibrate the model discrepancy parameters, we first observe that the dimension of $\t \in \R^{n_\theta}$ scales with the product of the state and decision variable dimensions. If the number of FOM evaluations is less than the decision variable dimension, i.e., $N < n_z$, then there exists infinitely many $\t$'s that fit the data~\eqref{eqn:discrepancy_approx}. Accordingly, we adopt a Bayesian calibration paradigm so that prior knowledge can constrain the calibration and uncertainty can be accounted for through the posterior. Bayesian calibration is challenging because $\d$ is
an operator mapping from the decision variable space to the state space.
We summarize the approach below and refer the reader to~\cite{hart2024hdsasampling,hart_vbw_2025} for additional
details.

\subsubsection*{Prior discrepancy}
The prior for $\t$ should be defined on the space of affine operators and
discretized to ensure that it respects the pre-discretization structures.
Following~\cite{hart2024hdsasampling,hart_vbw_2025}, we consider a
mean-zero Gaussian prior for the discrepancy parameters $\t$. Recall
that, for a general a mean-zero Gaussian random
vector $\vec{X}$, its probability density function is proportional to
$\exp \left( -\frac{1}{2} \| \vec{x} \|_{\W_{\vec{x}}} \right)$, where
$\W_{\vec{x}}$ is the precision matrix (the inverse of the covariance
matrix). We define the prior covariance for the discrepancy
parameters $\t$ by introducing a weighted norm $\| \cdot \|_{\W_\t}$
and take the weighting matrix $\W_\t$ as the prior precision.

To derive $\W_\t$, we define mean-zero Gaussian
priors on the state and decision variables spaces, resulting in
precision matrices $\W_\u$ and $\W_\z$ that define inner products for $\u$ and $\z$,
respectively. Considering an inner product on the space of affine operators,
we can derive an expression for $\W_\t$ that depends on $\W_\u$ and $\W_\z$.
These precision matrices are defined using theory and experience from
the analysis of infinite-dimensional Bayesian inverse problems~\cite{stuart_inv_prob,ghattas_infinite_dim_bayes_1}.
See~\cite{hart_vbw_2025} for detailed analysis of the prior specification.

\subsubsection*{Posterior discrepancy}
The data likelihood function is defined by assuming a mean-zero Gaussian noise
model. Since the data consists of full spatio-temporal observations of the
state discrepancy, we define the noise covariance as a scalar multiple
of the identity operator in the state space. This corresponds to a Gaussian noise model with covariance
$\alpha_d \M_\u^{-1}$, where $\alpha_d > 0$ scales the variance and
$\M_\u$ defines the state-space inner product\footnote{The
  matrix $\M_\u$ defines an inner product that integrates over both
  time and space. Hence, $\M_\u$ does not correspond to the mass
  matrix in the spatial discretization, but rather is a Kronecker
  product of a temporal weighting matrix and the spatial
  discretization mass matrix.}. In Bayesian inverse
problems, the noise level is determined based on characteristics of
the data-collection process.  In our context, the data is
  not noisy as is typical with sensor measurements.  Instead, data arises
  from simulations, and therefore the level of noise
  corresponds to linearization error from fitting an affine model to
  data from a nonlinear operator. We interpret $\alpha_d$
  as a linearization error magnitude and hence specify it relative to the
magnitude of the discrepancy data being fit.

Since we have assumed Gaussian prior and noise models, and parameterized $\d$ as a linear function of $\t$, the posterior distribution for $\t$ is Gaussian, with its mean and covariance matrix available in closed form. To simplify the exposition, we omit explicit expressions for the posterior mean and samples, but refer the reader
to~\cite{hart2024hdsasampling}, which demonstrates efficient posterior computation.

\subsection{Optimal solution updating}  \label{ssec:opt_sol_update}

Our goal is to approximate the~\texttt{FOMCO} solution online using the~\texttt{ROMCO} minimizer $\ztilde$
and $N=\mathcal O(1)$ FOM evaluations. To analyze how the model discrepancy
influences the optimization problem, observe that if $\t_0=\vec{0}$, then $\d(t;\z,\t_0) \equiv 0$ and the~\texttt{ROMCO} minimizer $\ztilde$ satisfies
 the first-order optimality condition
\begin{align}
\label{eqn:first_order_opt}
\nabla_\z \mathcal J(\ztilde,\t_0) = \vec{0},
\end{align}
where
\begin{align*}
\mathcal J(\z,\t) = J(\V \hat{\S}(t;\z) + \d(t;\z,\t),\z)
\end{align*}
 is the objective function of the~\texttt{ROMCO-MD} problem (which coincides with the~\texttt{ROMCO} objective function when $\t=\t_0$). We assume that $\mathcal J$ is twice continuously differentiable with respect to $(\z,\t)$ and that the Hessian $\nabla_{\z,\z} \mathcal J(\ztilde,\t_0)$ is positive definite. Applying the Implicit Function Theorem to~\eqref{eqn:first_order_opt} yields the existence of an operator $\z^\star: \mathcal N(\t_0) \to \R^{n_z}$, defined on a neighborhood of $\t_0$, such that
\begin{align*}
\nabla_\z \mathcal J(\z^\star(\t),\t) = \vec{0}
\end{align*}
for all $\t \in \mathcal N(\t_0)$. That is, the mapping $\z^\star(\t)$, from the model discrepancy parameters to the optimal decision variable, is uniquely defined in a neighborhood of the~\texttt{ROMCO} solution. Furthermore, the Jacobian of $\z^\star(\t)$, referred to as the post-optimality sensitivity operator, is given by
\begin{align*}
\nabla_\z \z^\star(\t_0) = - \H^{-1} \B
\end{align*}
where $\H = \nabla_{\z,\z} \mathcal J(\ztilde,\t_0)$ and $\B = \nabla_{\z,\t} \mathcal J(\ztilde,\t_0)$ are the $(\z,\z)$ and $(\z,\t)$ Hessians of $\mathcal J$, evaluated at $(\ztilde,\t_0)$. The linear approximation
\begin{align}
\label{eqn:linear_approx}
\z^\star(\t) \approx \z^\star(\t_0) + \nabla_\t \z^\star(\t_0) \t = \ztilde - \H^{-1} \B \t
\end{align}
provides a computationally efficient approach to update
the~\texttt{ROMCO} solution. Specifically, the decision variables can be updated
using $N$ FOM evaluations by
computing~\eqref{eqn:linear_approx}, where $\t=\overline{\t}$ is the posterior mean arising
from the Bayesian calibration. Uncertainty in this update can be quantified
by computing posterior samples for $\t$ and pushing them forward through~\eqref{eqn:linear_approx} to compute posterior decision variable samples.

\begin{figure}[h]
    \centering
    \includegraphics[width=0.8\textwidth]{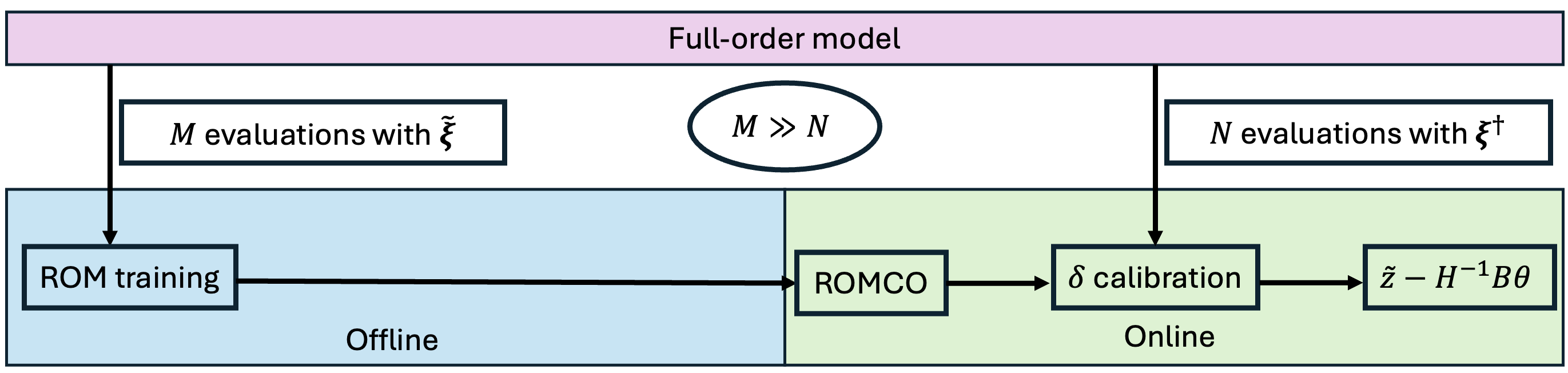}
    \caption{Summary of the HDSA-MD approach for updating~\texttt{ROMCO} solutions.}
    \label{fig:summary}
\end{figure}

We summarize our HDSA-MD approach in~\Cref{fig:summary}.

\section{Atmospheric contaminant control using operator inference} \label{sec:atmosphere_transport}

  This section presents a control
problem for a two-species advection-diffusion-reaction system.
Let $u_1(\x,t)$ denote a contaminant that is to be neutralized by deploying a decontaminant $u_2(\x,t)$ within a spatial domain $\Omega$.
Both $u_1$ and $u_2$ diffuse homogeneously and advect along a stationary and spatially heterogeneous velocity field $\v$, which is uncertain (later parameterized by $\xivec$).
The following system of PDEs describes this scenario:

\begin{subequations}
\label{eq:adr-pde}
\begin{align}
    & \frac{\partial u_1}{\partial t}
    = \kappa_1\Delta u_1 - \v\cdot \nabla u_1 - \rho u_1 u_2 & \qquad \text{on } \Omega \times (0,T],
    \\
    & \frac{\partial u_2}{\partial t}
    = \kappa_2\Delta u_2 - \v\cdot \nabla u_2 - \rho u_1 u_2 + s & \qquad \text{on } \Omega \times (0,T],
    \\
    \label{eq:adr-bcs}
    & \vec{n}\cdot\nabla u_1 = \vec{n}\cdot\nabla u_2 = 0 & \qquad \text{on } \partial \Omega \times (0,T],
    \\
   &  u_1(\x, 0) = \mathring{u}_1(\x),
    \quad
    u_2(\x, 0) = 0 & \qquad \text{on } \Omega .
\end{align}
\end{subequations}

Here, $T=0.4$ is the final time, $\kappa_1=\kappa_2 =0.1$ are diffusion coefficients, and $\mathring{u}_1$ is the initial contaminant profile.
The homogeneous Neumann boundary
conditions~\cref{eq:adr-bcs} prescribe zero-flux
boundaries. The decontaminant enters the domain through the source term
$s$ and reacts with the contaminant with rate $\rho=2$.

\begin{figure}[h]
\centering
\begin{subfigure}
    \centering
    \includegraphics[height=.225\textheight]{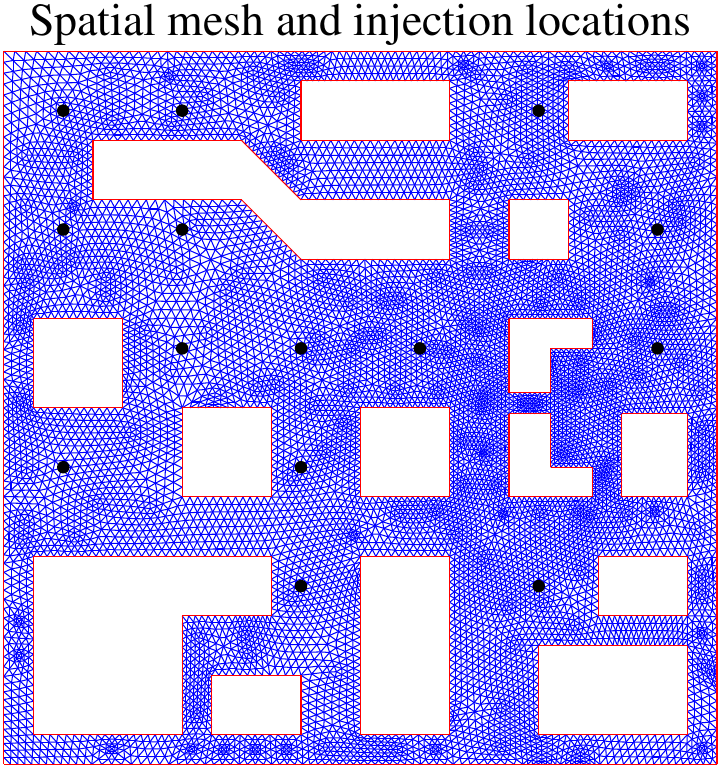}
\end{subfigure}
\hspace{.1cm}
\begin{subfigure}
    \centering
    \includegraphics[height=.225\textheight]{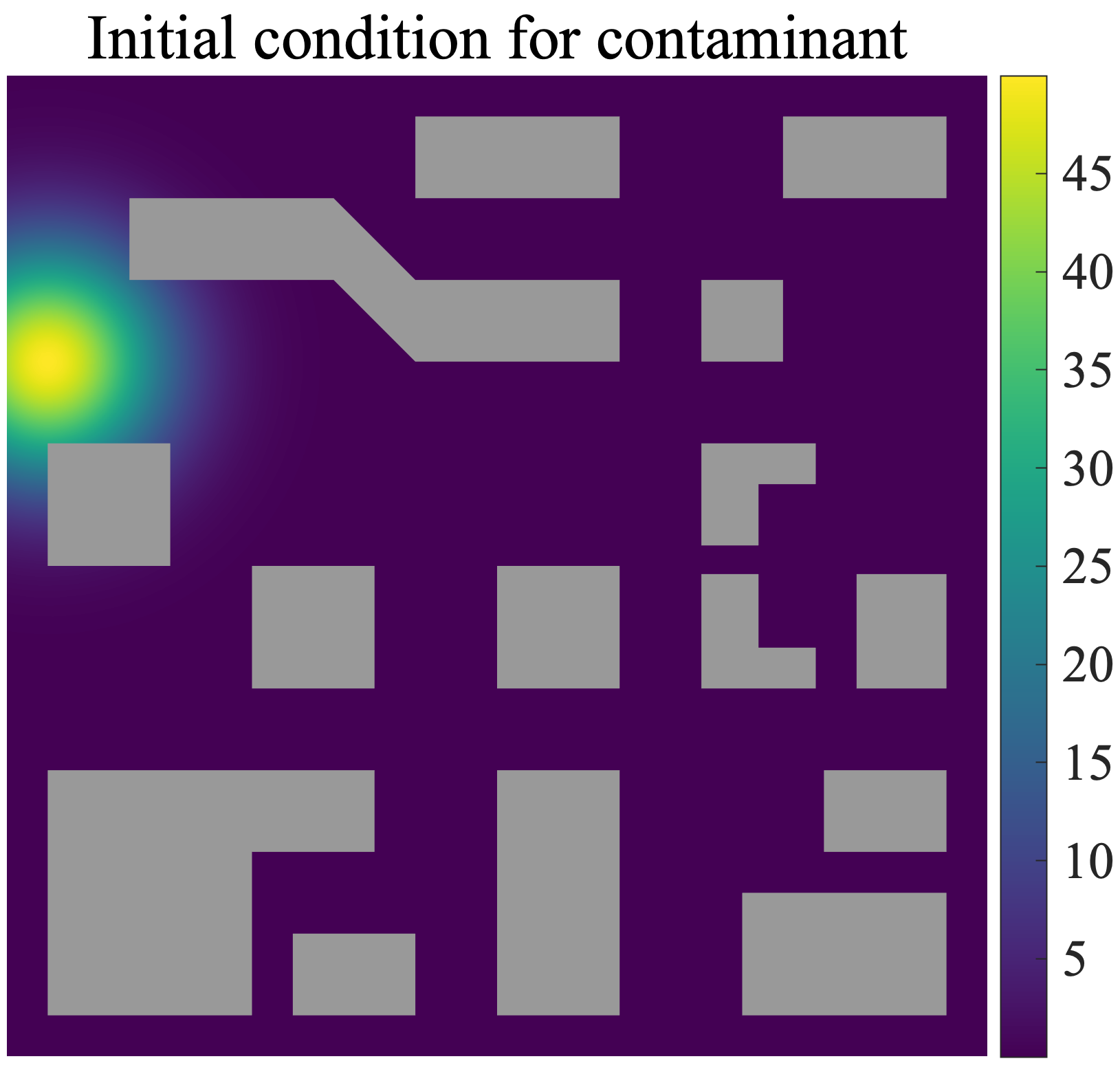}
\end{subfigure}
\hspace{.1cm}
\begin{subfigure}
    \centering
    \includegraphics[height=.225\textheight]{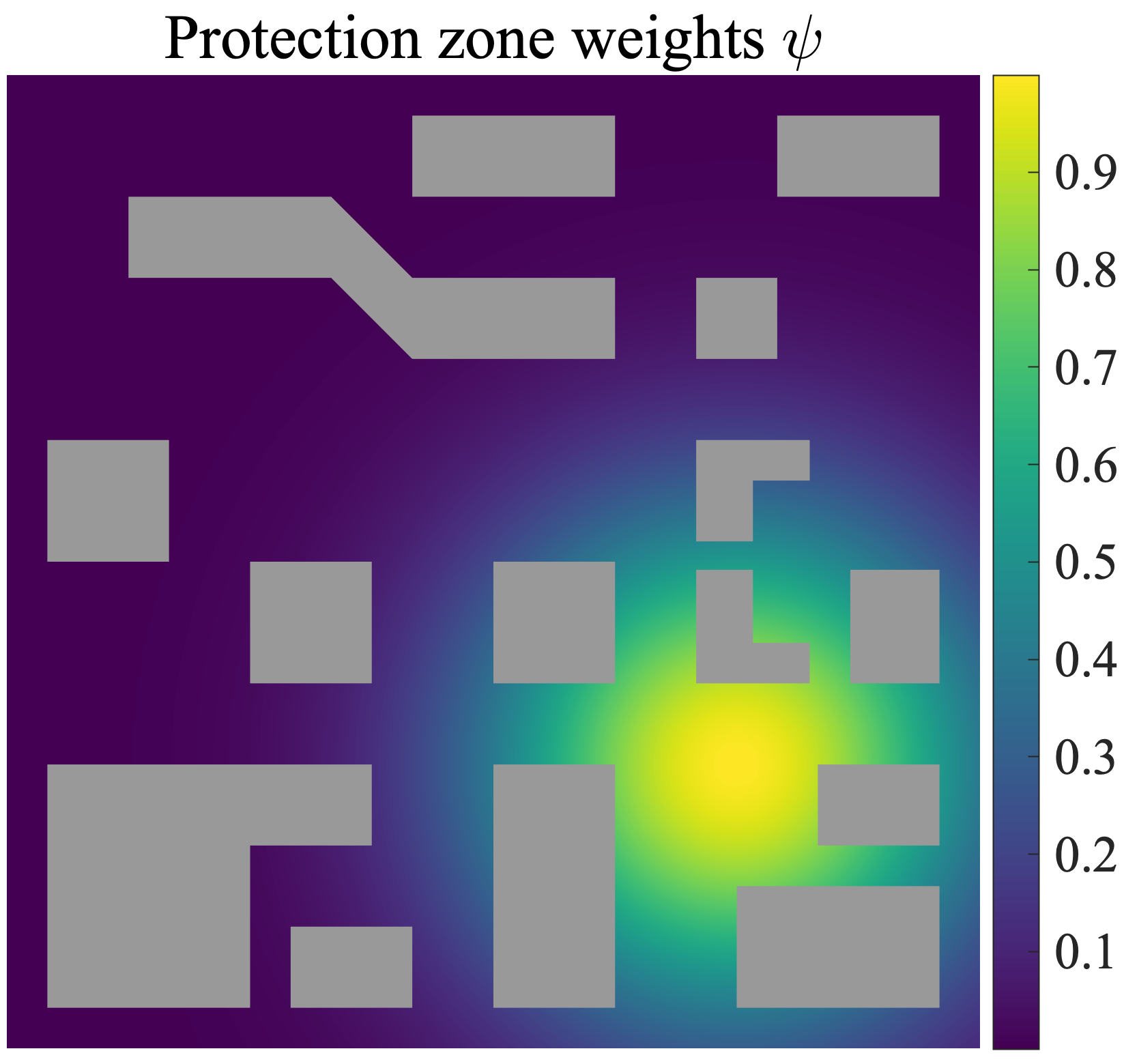}
\end{subfigure}
\caption{Left: finite element mesh for the spatial domain $\Omega$. The decontaminant injection points $\bar{\x}_1,\ldots,\bar{\x}_{n_q}$ are marked in black. Center: initial condition $\mathring{u}_1$ for the contaminant. Right: weight function $\psi(\x)$ for measuring the amount of contaminant near the protection zone.
}
\label{fig:adr-mesh}
\end{figure}

We consider a two-dimensional domain $\Omega \subset (0,1.2)^2$ that excludes obstacles and discretize in space with finite
elements using $n_x = 10{,}702$ nodes. This leads to a system of
ordinary differential equations
\begin{subequations}
\label{eq:adr-ode}
\begin{align}
    \M\dot{\u}_1(t)
    &= \A_1\u_1(t) + \mathbf{R}_1\left( \u_1(t) \otimes \u_2(t) \right),
    \label{eq:adr-ode-u1}
    \\
    \M\dot{\u}_2(t)
    &= \A_2\u_2(t) + \mathbf{R}_2 \left( \u_1(t) \otimes \u_2(t) \right) + \M\s(t),
    \label{eq:adr-ode-u2}
    \\
    \u_1(0) &= \mathring{\u}_1,
    \quad
    \u_2(0) = \vec{0},
\end{align}
\end{subequations}
where $\u_1(t),\u_2(t),\s(t),\mathring{\u}_1 \in\R^{n_x}$ are,
respectively, the spatial discretizations of $u_1$, $u_2$, $s$, and
$\mathring{u}_1$; $\M\in\R^{n_x\times n_x}$ is the finite-element mass
matrix. The matrices $\A_1,\A_2 \in \R^{n_x \times n_x}$ and
$\mathbf{R}_1,\mathbf{R}_2 \in \R^{n_x \times n_x n_x}$ correspond to the
discretization of the linear and bilinear operators
in~\eqref{eq:adr-pde}, respectively. The spatial domain, mesh, and initial state are shown in \Cref{fig:adr-mesh}.

\subsection{Control problem formulation}

The goal is to manipulate the source term $\s(t)$ to minimize the
contaminant $u_1$ in a
protection region that is represented by $\psi(\x) = \exp\left(-10\|\x
  - \x_p\|^2\right)$, where $\x_p = (0.9, 0.35)$. We seek to minimize
 the contaminant concentration in the protection
region, measured as $\| u_1(\x,t) \psi(\x) \|$, by injecting a decontaminant at $n_q = 14$ discrete spatial points $\bar{\x}_1,\ldots,\bar{\x}_{n_q}\in\Omega$.
The protection region and injection locations are shown in~\Cref{fig:adr-mesh}.

We model the source term as
a sum of small-variance Gaussians centered at the injection points $\bar{\x}_1,\ldots,\bar{\x}_{n_q}$,
\begin{align}
\label{eqn:decontaminant_source}
    s(\x,t)
    = \sum_{i=1}^{n_q}q_{i}(t;\z)^2 \phi_i(\x),
    \qquad\qquad
    \phi_i(\x) = 50 \exp\left(-1000\|\x - \bar{\x}_i\|^2\right),
\end{align}
where $q_i(t;\z)$, $i=1,2,\dots,n_q$, are functions of time parameterized by the decision variable $\z$.
Note that the source magnitude is defined by $q_i(t;\z)$ squared to ensure that the injection is
nonnegative.

Letting $\u(t) = (\u_1(t)\trp~\u_2(t)\trp)\trp\in\R^{2n_x}$, $\vec{\psi} \in \R^{n_x}$, and $\q(t;\z)= (q_1(t;\z)\trp,q_2(t;\z)\trp,\dots,q_{n_q}(t;\z)\trp)\trp\in\R^{n_q}$ denote the discretized state, protection zone weights, and injection magnitudes, respectively, we consider the optimization problem
\begin{subequations}
\label{eq:adr-opt}
\begin{gather}
    \label{eq:adr-objective}
    \min_{\z \in \R^{n_z}}
    \frac{1}{2}\int_{0}^{T}\|\u_1(t) \odot \vec{\psi} \|_{\M}^2
    + \gamma\|\q(t;\z)\odot\q(t;\z)\|^2\,dt
    \\
    \label{eq:adr-constraint}
    \textrm{where}~~(\u(t),\q(t;\z))~~\textrm{ jointly solves~\cref{eq:adr-ode}},
\end{gather}
\end{subequations}
where $\odot$ denotes the Hadamard (elementwise) product and $\z \in \R^{n_z}$ defines the $n_q$ injection magnitudes at a collection of $n_s$ time nodes, $n_z=n_q n_s$. The regularization term $\gamma \|\q(t;\z)\odot\q(t;\z)\|^2 = \gamma \sum_{i=1}^{n_q} q_i(t;\z)^4$ induces a cost, scaled by the regularization coefficient $\gamma = 10^{-5}$, to inject the decontaminant.\footnote{Note that the regularization involves the norm of $\q(t;\z)\odot\q(t;\z)$ since regularizing with the norm of $\q(t;\z)$ would be equivalent to an optimization problem involving $\ell^1$ regularization and non-negativity constraints, thus leading to a lack of differentiability that hinders the post-optimality analysis.}

\begin{figure}[t]
    \centering
    \includegraphics[width=0.4\textwidth]{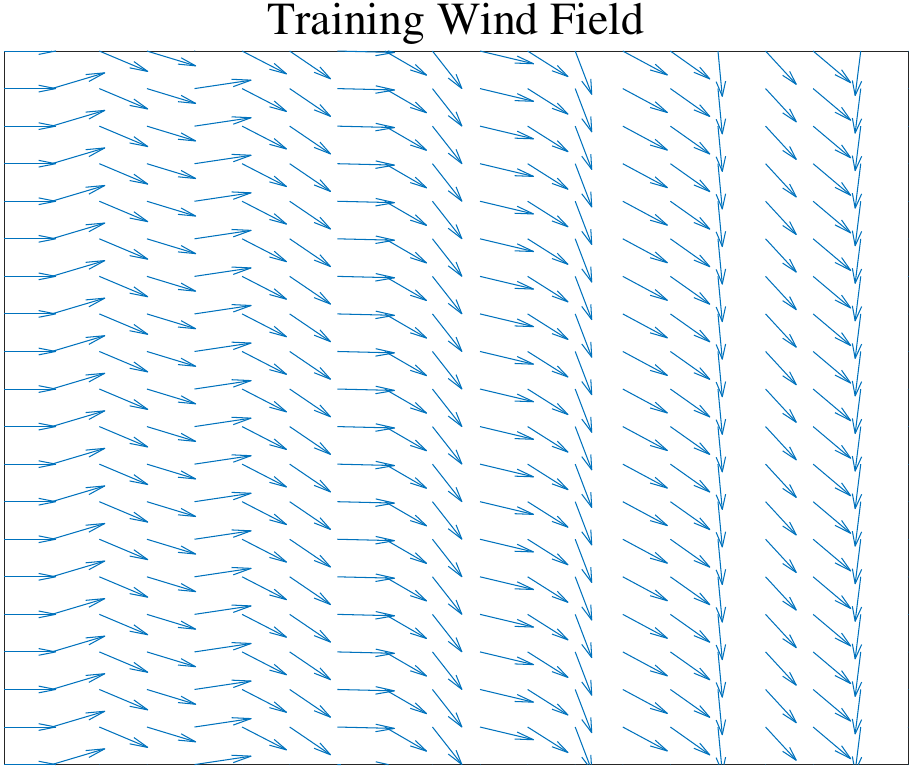}
    \hspace{.5cm}
    \includegraphics[width=0.4\textwidth]{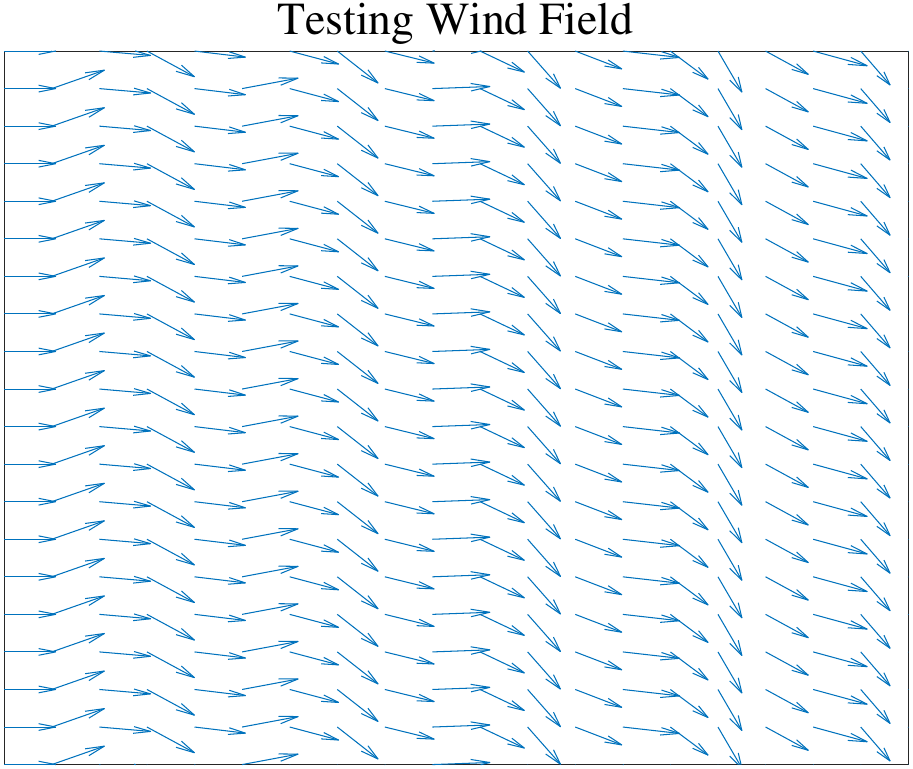}
    \caption{Wind fields used for training, $\tilde{\xivec} = (1,1,10,75)\trp$, and testing, $\xivec^\dagger = (1.5,1.5,8,50)\trp$.}
    \label{fig:adr-winds}
\end{figure}

Our goal is to compute (or approximate) the solution of~\eqref{eq:adr-opt} in real-time. However, the ODE system~\eqref{eq:adr-ode} depends on the velocity field $\v$, which is only known online when the optimization solution is required. For this experiment, the velocity field is parameterized by the vector $\xivec = (\xi_1,\xi_2,\xi_3,\xi_4) \in \R^4$ as
\begin{align}
    \label{eq:adr-velocity}
    \v(\x; \xivec)
    = \left[\begin{array}{c}
        4\xi_1\cos(\xi_4 x_1/L) - 2\xi_2\sin(\xi_3 \pi x_1)\sin(\xi_4 x_1/L) \\
        4\xi_1\sin(\xi_4 x_1/L) - 2\xi_2\sin(\xi_3 \pi x_1)\cos(\xi_4 x_1/L)
    \end{array}\right].
\end{align}
We assume an offline training wind field parameter vector $\tilde{\xivec} = (1,1,10,75)\trp$. During the online phase, the testing wind field parameter vector is $\xivec^\dagger = (1.5,1.5,8,50)\trp$. The two wind fields are shown in \Cref{fig:adr-winds}.

We will train an operator inference ROM offline using the training wind field parameter vector $\tilde{\xivec}$. In the online phase, we solve the optimal control problem using the ROM, compute a single full-model evaluation~\eqref{eq:adr-ode} with testing wind field parameter vector $\xivec^\dagger$, and update the optimal control using HDSA-MD. The subsections that follow detail these steps.

\subsection{Operator inference ROM}
\label{sec:adr-opinfrom}

Training data is generated by constructing $n_c = 5$ smooth (as functions of time) random
control vectors $\q^{(\ell)}(t)$, $\ell=1,\ldots,n_c$, and solving the
full-order system~\cref{eq:adr-ode} with these controls and the
training wind field parameter vector $\tilde{\xivec}$. The training
controls are shown in the left panel of
\Cref{fig:adr-controls}.  A proper orthogonal decomposition (POD)~\cite{sirovich1987pod} is computed for the two states separately to produce basis matrices
$\V_1 \in \R^{n_x \times r_1}$ and $\V_2 \in \R^{n_x \times r_2}$ so
that the states are approximated as
\begin{align}
    \label{eq:adr-approx}
    \u_1(t) \approx \V_1\widehat{\u}_1(t),
    \qquad
    \u_2(t) \approx \V_2\widehat{\u}_2(t),
\end{align}
for $\widehat{\u}_1(t)\in\R^{r_1}$ and $\widehat{\u}_2(t)\in\R^{r_2}$. In our example, $r_1=19$ and $r_2=33$ are determined as the smallest basis dimensions such that the POD residual energy is less than $10^{-5}$; see~\eqref{eq:adr-podres} for details.

Inserting~\cref{eq:adr-approx} into~\cref{eq:adr-ode} and left multiplying the $k$th equation by $\V_k\trp$, we obtain a system of ODEs
\begin{subequations}
\label{eq:adr-rom}
\begin{align}
    \dot{\widehat{\u}}_1(t)
    &= \widehat{\A}_1\widehat{\u}_1(t) + \widehat{\mathbf{R}}_1\left( \widehat{\u}_1(t)\otimes\widehat{\u}_2(t) \right),
    \\
    \dot{\widehat{\u}}_2(t)
    &= \widehat{\A}_2\widehat{\u}_2(t) + \widehat{\mathbf{R}}_2 \left( \widehat{\u}_1(t)\otimes\widehat{\u}_2(t) \right) + \widehat{\vec{\Phi}} \left( \q(t) \odot \q(t) \right),
\end{align}
\end{subequations}
where the system matrices are given by
\begin{align}
    \label{eq:adr-galerkin-matrices}
    \begin{aligned}
    \widehat{\A}_k &= \V_k\trp\A_k\V_k \in\R^{r_k\times r_k},
    \\
    \widehat{\mathbf{R}}_k &= \V_k\trp\mathbf{R}_k(\V_1\otimes\V_2)\in\R^{r_k \times r_1 r_2},
    \end{aligned}
    \qquad k=1,2,~\text{and}\qquad
    \widehat{\vec{\Phi}} = \V\trp\M \vec{\Phi} \in\R^{r_2 \times n_q}.
\end{align}
Here, $\vec{\Phi} \in \R^{n_x \times n_q}$ corresponds to evaluating the injection basis functions $\phi_i$ from~\eqref{eqn:decontaminant_source}, $i=1,2,\dots,n_q$, at the spatial nodes.

Classical model-reduction methods use the matrices of the FOM~\cref{eq:adr-ode} to construct the ROM matrices~\cref{eq:adr-galerkin-matrices}, but this requires intrusive access to the FOM code to explicitly construct $\M,\A_1,\A_2,\mathbf{R}_1,\mathbf{R}_2,$ and $\vec{\Phi}$. We construct a ROM using operator inference~\cite{ghattas2021physicsbased,kramer2024opinfsurvey,peherstorfer2016opinf}, a non-intrusive alternative to classical methods. Specifically, operator inference postulates the model structure~\eqref{eq:adr-rom} and estimates the matrices~\eqref{eq:adr-galerkin-matrices} via least-squares regression using the training data snapshots. Details of the regression problem formulation and solution are given in Appendix~\ref{appendix:atmosphere_transport}.

\subsection{ROM-constrained optimization}

Given our operator inference ROM, we consider the \texttt{ROMCO} problem,
\begin{subequations}
\label{eq:adr-opt-rom}
\begin{gather}
    \label{eq:adr-objective-rom}
    \min_{\z \in \R^{n_z}}
    \frac{1}{2}\int_{0}^{T}\|\V_1\widehat{\u}_1(t) \odot \vec{\psi} \|_{\M}^2
    + \gamma\|\q(t;\z)\odot\q(t;\z)\|^2\,dt
    \\
    \label{eq:adr-constraint-rom}
    \textrm{where}~~(\widehat{\u}(t),\q(t;\z))~~\textrm{ jointly solves~\cref{eq:adr-rom}},
\end{gather}
\end{subequations}
where $\widehat{\u}(t) = (\widehat{\u}_1(t)\trp~~\widehat{\u}_2(t)\trp)\trp\in\R^{r_1 + r_2}$.

\begin{figure}[t]
    \centering
    \begin{subfigure}
        \centering
        \includegraphics[width=.45\textwidth]{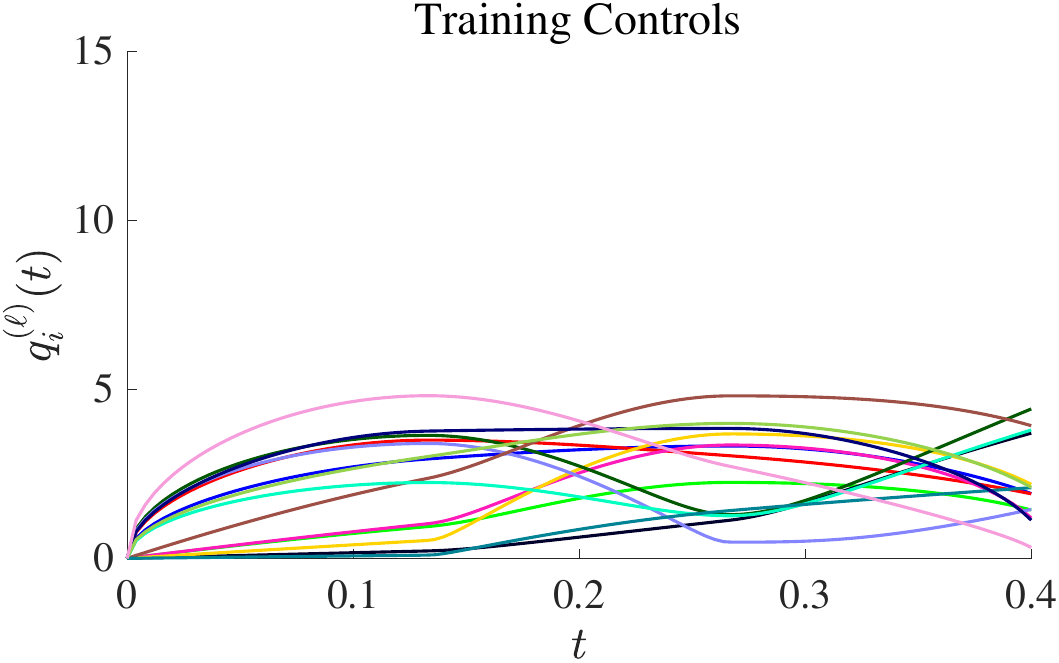}
    \end{subfigure}
    \hfill
    \begin{subfigure}
        \centering
        \includegraphics[width=.45\textwidth]{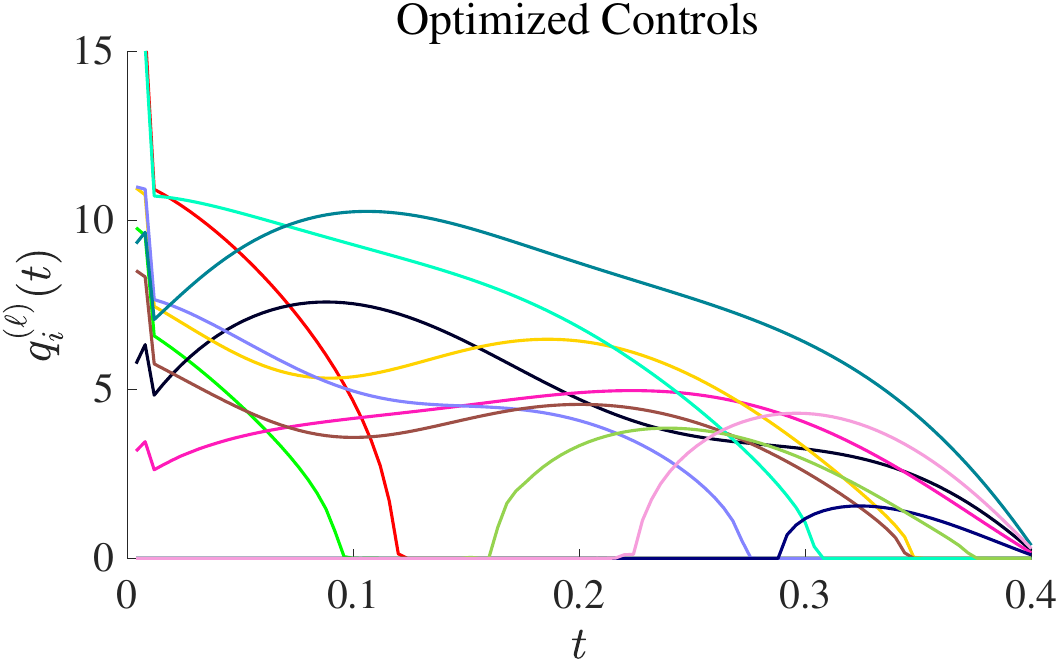}
    \end{subfigure}
    \caption{Left: components of one of the random training controls $\q^{(\ell)}(t)$ used to generate training data.
    Right: components of the optimized controls, the solution of the \texttt{ROMCO} problem~\cref{eq:adr-opt-rom}.
    }
    \label{fig:adr-controls}
\end{figure}

The discretized control is defined at $n_s=100$ time nodes and $n_q=14$ spatial locations, hence the optimization algorithm must search a $n_z=1400$ dimensional space. The \texttt{ROMCO} problem~\cref{eq:adr-opt-rom} is solved using derivative-based optimization, where objective-function gradient evaluations and Hessian-vector products are computed using the adjoint method. This enables efficient computation of derivatives in that the number of ROM evaluations per derivative evaluation is independent of the decision variable dimension $n_z$. The right panel of \Cref{fig:adr-controls}~shows the optimized controls $\q(t;\z)$. Observe that some components are nonzero over the entire time interval while others are only nonzero over a subinterval.

\begin{figure}[t]
    \centering
    \includegraphics[width=.9\textwidth]{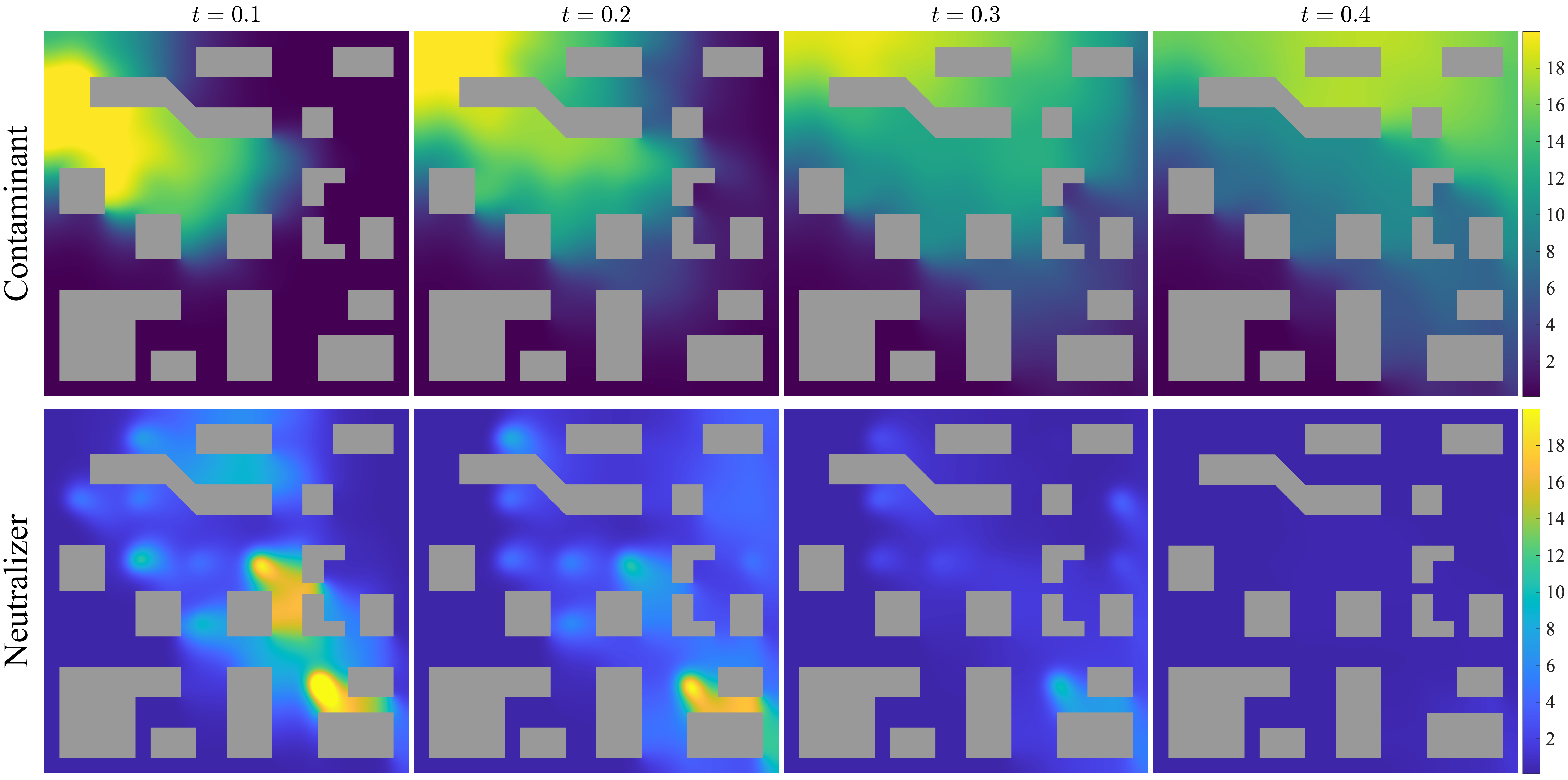}
    \caption{Solution to the full-order model with the wind field parameter vector $\xivec=\xivec^\dagger$ and the control obtained from solving the \texttt{ROMCO} problem~\cref{eq:adr-opt-rom}.}
    \label{fig:adr-opt}
\end{figure}

The controls shown in~\Cref{fig:adr-controls} are optimal for the ROM constructed using the training wind field defined by $\xivec = \tilde{\xivec}$. We seek a solution in the test case when the wind field is defined by $\xivec=\xivec^\dagger$.  To assess the performance of the \texttt{ROMCO} problem solution, we simulate the system using the full-order model with the testing wind field and \texttt{ROMCO} controls. The state solution is depicted in \Cref{fig:adr-opt}, where we observe that a nontrivial amount of contaminant reaches the protection zone.

\subsection{Model discrepancy analysis}

As seen in \Cref{fig:adr-opt}, the \texttt{ROMCO} control solution is insufficient. We seek to improve the controls using only a single evaluation of the FOM solution operator $\S(t;\tilde{\z},\xivec^\dagger)$. As outlined in~\Cref{ssec:dis_calibration}, we use this single FOM evaluation to calibrate a representation of the discrepancy $\d(t;\z,\t) \approx \S(t;\z,\xivec^\dagger)-\V \widehat{\S}(t; \z)$, where $\V$ is the concatenation of the basis matrices $\V_1$ and $\V_2$, and $\widehat{\S}(t; \z)$ is the solution operator for the ROM~\eqref{eq:adr-rom}.

The decision-variable dimension is $n_z=1400$ and yet we have only $N=1$ FOM evaluation. Consequently, there is significant posterior uncertainty in the model discrepancy
operator. However, notice that the objective function $\mathcal J(\z,\t_0)$
has small curvature in many directions. Consider a Taylor approximation
\begin{align*}
\mathcal J(\z,\t_0) \approx \mathcal J(\ztilde,\t_0) + \frac{1}{2} (\z-\ztilde)^T \H (\z-\ztilde)
\end{align*}
where the linear term is omitted since $\nabla_\z \mathcal J(\ztilde,\t_0)=\vec{0}$ and the higher-order (greater than quadratic) terms are neglected.

Let $\H \vec{w}_i=\lambda \vec{w}_i$ denote the eigenvalues and eigenvectors of the Hessian $\H$, $i=1,2,\dots,n_z$. We assume the eigenvalues are ordered as $\lambda_1 \ge \lambda_2 \ge \cdots \ge \lambda_{n_z}$. For a small eigenvalue $\lambda_i$, directions for which $\z-\ztilde$ align with the eigenvector $\vec{w}_i$ will have a negligible $\mathcal O(\lambda_i)$ magnitude influence on the objective function. However, observing~\eqref{eqn:linear_approx}, we see that $\t$'s such that $\B \t$ aligns with $\vec{w}_i$ will result in optimal solution updates that scale like $\lambda^{-1}$.  Hence, if posterior samples align with eigenvectors corresponding to small eigenvalues of the Hessian, then we will observe large variances in optimal solution updates despite the fact that these large variances yield small changes in the objective function. To avoid such scenarios, we consider an alternative form of the optimal solution update that projects the post-optimality sensitivity onto the subspace of the Hessian's leading eigenvectors.

Let $\P \in \R^{n_z \times n_z}$ be the orthogonal projector onto the span of the $r$ leading eigenvectors $\{ \vec{w}_j \}_{j=1}^r$. The projected optimal solution update is given by
\begin{align*}
\ztilde - \P \H^{-1} \B \t .
\end{align*}
In this example, we retain the $r=15$ leading eigenvectors of the Hessian as they exhibit two orders of magnitude decay, i.e., $\frac{\lambda_{15}}{\lambda_1} =\mathcal O(10^{-2})$.

\begin{figure}[t]
    \centering
    \includegraphics[width=0.49\textwidth]{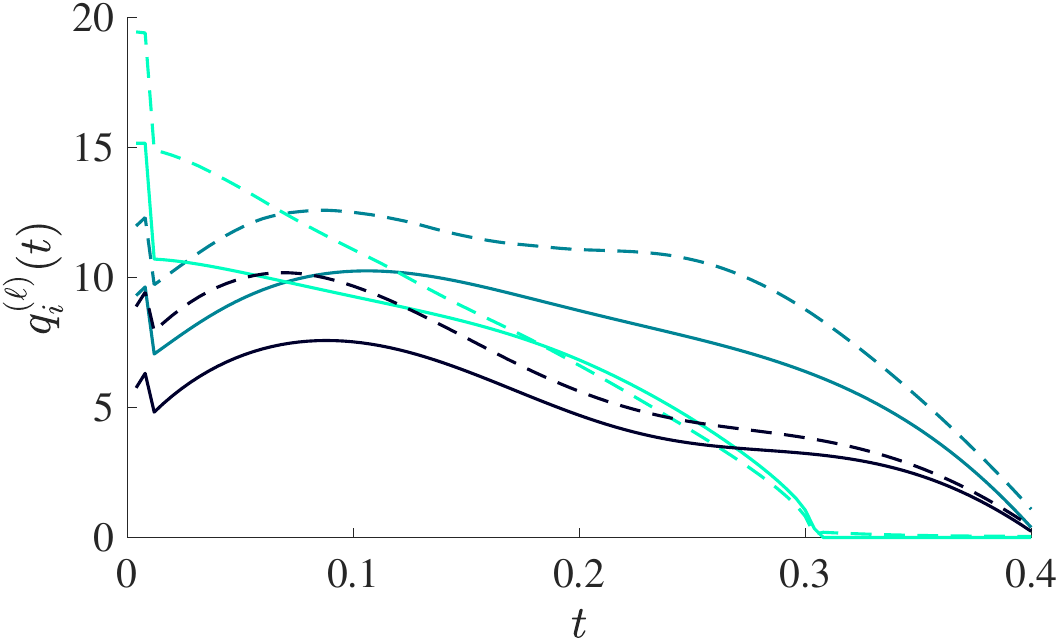}
      \includegraphics[width=0.49\textwidth]{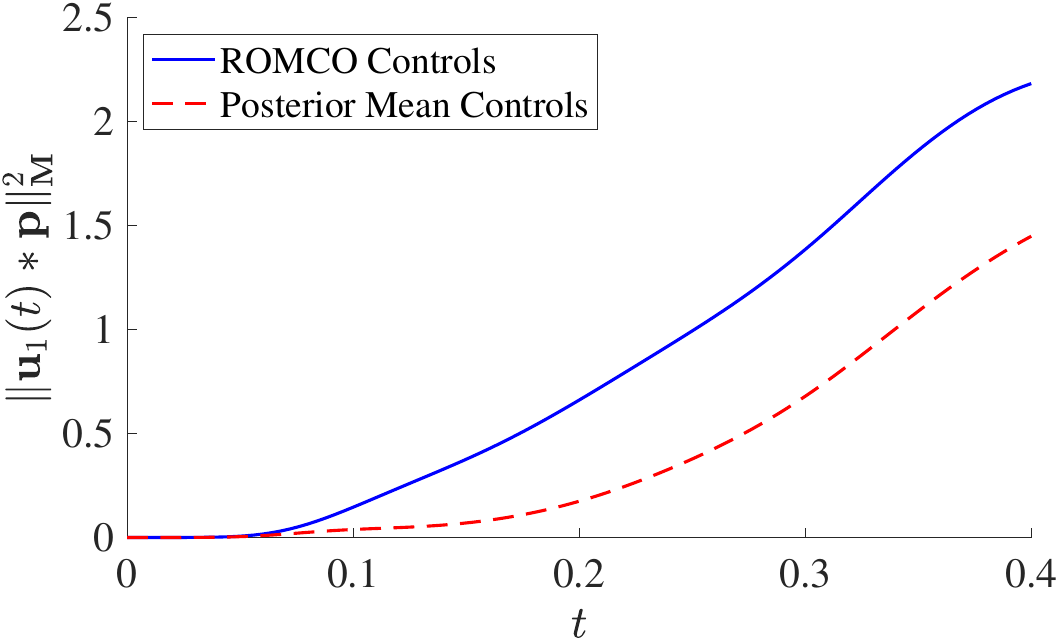}
    \caption{Left: three components of the optimized controls, the solution of the \texttt{ROMCO} problem~\cref{eq:adr-opt-rom} (in solid lines) and the posterior mean controls (in broken lines). Right: time evolution of the contaminant mass in the protection region using the optimized controls and the posterior mean controls.}
    \label{fig:adr-update}
\end{figure}

The left panel of \Cref{fig:adr-update} displays three components of the posterior mean controls alongside the optimized controls that were computed from the \texttt{ROMCO} problem~\cref{eq:adr-opt-rom}. These three components were chosen because they had the largest change (in the post-optimality update) amongst the $n_q=14$ components. The right panel of \Cref{fig:adr-update} shows the time evolution of the protection region contaminant mass squared, $\|\u_1(t)\odot\vec{\psi}\|_{\M}^2$, which is the integrand in \eqref{eq:adr-objective}. The reduction in contaminant demonstrates the benefits gained from using HDSA-MD with only a single full-order-model evaluation. Quantitatively, we can measure the control improvement by comparing the values of the \texttt{FOMCO} objective function. This requires an addition FOM evaluation that is typically not computed in practice, but done here for the purpose of validation. The \texttt{FOMCO} objective function, evaluated with the wind field parameter vector $\xivec^\dagger$ and controls defined by $\tilde{\z}$ and, is $0.355$. In comparison, the \texttt{FOMCO} objective function is $0.233$ when evaluated with the wind field parameter vector $\xivec^\dagger$ and controls defined by $\tilde{\z}-\P \H^{-1} \B \overline{\t}$, where $\overline{\t}$ denotes the posterior mean of the model discrepancy parameters. Hence, we realized a $34\%$ reduction in the \texttt{FOMCO} objective using only a single FOM evaluation to update the controls.

\section{Wildfire ignition inversion using a neural operator ROM} \label{sec:wildfire}

This section presents an inverse problem for the ignition point of a wildfire in the Valles Caldera of northern New Mexico, which features heterogeneous topography and fuels. Additionally, the caldera is a high-consequence area due to nearby infrastructure, its status as a protected nature preserve, and a history of wildfire activity. Such inversion supports initialization of fire propagation codes as well as forensic analysis~\cite{shaddy2024generative}. Wildfire modeling is complex as it involves a system of nonlinear PDEs that couples the atmosphere and fire propagation. We use the Fortran package WRF-SFIRE~\cite{mandel2011coupled} to simulate the fire dynamics, as summarized below.

Let $u(\x,t)$ denote the {\em level-set function}, defined on a spatial domain $\Omega \subset \R^2$. Physical quantities resulting from the fire---such as temperature, heat flux, and the fire front---are computed from $u$ in a postprocessing phase. The interpretation of the level-set function is a signed distance to the fire front. Specifically, at a given time $t$, the points $\x \in \Omega$ where $u(\x,t)=0$ define the fire front, the points where $u(\x,t)<0$ define the burned region, and the points where $u(\x,t)>0$ define the unburned region. The magnitude of the level-set function measures distance from the fire front.

The wildfire is modeled using the coupled PDEs
\begin{subequations}
\label{equ:wrf-dynamics}
\begin{alignat}{2}
	& \frac{\partial u}{\partial t} + \Psi(\vec{v}) \| \nabla u\| = 0, \qquad && \text{on }  \Omega \times (0,T], \label{equ:wrf-lfn} \\
	& \frac{D \vec{v}}{D t} = f(u), && \text{on }  \Omega \times (0,C) \times  (0,T], \label{equ:wrf-euler}
\end{alignat}
\end{subequations}
where $\vec{v}$ is three-dimensional atmospheric wind, $C>0$ is the model top of atmosphere,
$\frac{D}{Dt}$ is the material derivative, and $T>0$ is the final time. The fire spread rate $\Psi$
depends on $\vec{v}$ as well as material properties such as fuel type, topographic slope, and fuel moisture. Equation~\cref{equ:wrf-lfn} is the level-set equation, while~\eqref{equ:wrf-euler} are the Euler equations for atmospheric wind $\vec{v}$. There is a two-way coupling between the fire and atmosphere. Winds will drive the fire in a particular direction via the dependence of $\Psi$ on $\vec{v}$, while the heat from the fire will create updraft that affects winds via the Euler equation forcing $f(u)$.

To find a solution to~\cref{equ:wrf-dynamics}, we must specify initial conditions and, as appropriate, boundary conditions. The initial condition for the level-set function $u$ is defined by the shifted pointwise Euclidean distance to the ignition point,
\begin{align} \label{eqn:level_set_init}
\mathring{u}(\vec{x}; \vec{z}) = \| \vec{x} - \vec{z} \| + c,
\end{align}
where $\z \in \R^2$ is the fire ignition location and $c=-10$ is a
shifting factor determined by the initial burn radius. Since $u$ is a
signed distance function, no boundary conditions are imposed on the
level-set function~\cite{mandel2011coupled}. The initial condition for the wind velocity is user-specified based on knowledge of the atmosphere, and open boundary conditions are employed; see~\cite{skamarock2019wrf} for details.

The Fortran package WRF-SFIRE~\cite{mandel2011coupled} handles all  physical properties, empirical relationships, and numerics for~\cref{equ:wrf-dynamics}. The level-set function $u(\x,t)$ is discretized using a $173 \times 133$ computational grid with spatial cells of size $60 \times 60$ meters, resulting in $n_u=23009$ spatial nodes. We consider a time horizon of $T=8$ hours with a time-step of $2$ seconds. WRF-SFIRE employs a finite-volume spatial discretization with the $(x,y)$ mesh for~\cref{equ:wrf-euler} being five times coarser than~\cref{equ:wrf-lfn}. Runge--Kutta time-steppers of order 2 and 3 are used for~\cref{equ:wrf-lfn} and~\cref{equ:wrf-euler}, respectively. Operator splitting couples the equations, where $u$ is updated first and then provided as an input to $\vec{v}$.

\subsection{Inverse problem formulation}

Our goal is to estimate the ignition location $\z \in \R^2$ of the fire, given observations of the burn area at a limited number of time points $\{t_j^{\text{obs}}\}_{j=1}^{n_\text{obs}}$ that occur after ignition. We assume hourly observations that begin two hours after the fire ignition, i.e., observations at times $t_j^{\text{obs}}=1+j$ hours, $j=1,\dots,n_\text{obs}$, where $n_\text{obs}=7$.

In practice, satellite observations enable assessment of whether a spatial location has been burned, i.e., binary data indicating that a pixel is burned or unburned. Such data may be post-processed to estimate the fire front and compute the level-set function from it. Accordingly, we consider observational data of the form $\vec{Y} = (\y_1, \dots, \y_{n_\text{obs}}) \in \R^{n_u \times n_\text{obs} }$, where $Y_{ij}$ is the observed signed-distance function at the $i^{th}$ spatial node $\x^{(i)}$ and time $t_j^{\text{obs}}$. Specifically, given the fire front $\Omega_f(t) \subset \R^2$ at time $t$, we denote the distance from a point $\vec{x} \in \R^2$ to the fire front as
$$
d(\vec{x}, \Omega_f(t)) = \inf_{\vec{y} \in \Omega_f(t)} \| \vec{x}-\vec{y} \| \,
$$
and define
\begin{equation}
Y_{ij} = \begin{cases}
    \displaystyle -d(\vec{x}^{(i)}, \Omega_f(t_j^{\text{obs}})) &\qquad \text{node}~\vec{x}^{(i)}~\text{is burned at}~t_j^{\text{obs}} \\
    \displaystyle ~~d(\vec{x}^{(i)}, \Omega_f(t_j^{\text{obs}})) &\qquad\text{otherwise}
\end{cases} \ .
\label{equ:wrf-obs}
\end{equation}

The ignition point is estimated by solving the optimization problem
\begin{equation}
	 \min_{\z \in \R^2} \sum_{j=1}^{n_\text{obs}} \| \u(t_j^{\text{obs}}; \z) - \y_j \|^2
	\label{equ:wrf-objective}
\end{equation}
where $\u(t_j^{\text{obs}};\z) \in \R^{n_u}$ denotes the solution of~\eqref{equ:wrf-lfn} after spatial discretization, evaluated at time $t=t_j^\text{obs}$, with the level-set initial condition parameterized by $\z$ as defined in~\eqref{eqn:level_set_init}.

Solving~\eqref{equ:wrf-objective} in real-time is intractable because the computational expense of WRF-SFIRE simulations limits the number of evaluations. Furthermore, it is impractical to implement sensitivity or adjoint equations due to the complex nonlinear feedbacks, operator splitting, and fire spread model characteristics used in WRF-SFIRE. Constructing an intrusive ROM would be challenging for the same reasons. Accordingly, we seek to evaluate WRF-SFIRE offline and train a non-intrusive ROM to approximate the level-set evolution.

\subsection{Neural operator ROM}
For the atmospheric transport problem in~\Cref{sec:atmosphere_transport}, we were able efficiently train a nonintrusive ROM via
operator inference by leveraging the quadratic structure of the FOM. Because WRF-SFIRE involves more complex nonlinear
structure and feedback mechanisms between the fire spread and atmospheric models, in this experiment we construct a nonintrusive
ROM through neural operators.

Neural operators are a large class of methods that seek to leverage
the strengths of neural networks to approximate mappings between
infinite-dimensional spaces~\cite{kovachki2023neural}. Many
contemporary methods fall into this category, including deep
operator networks (DeepONet)~\cite{lu2021learning}, Fourier neural
operators (FNO)~\cite{li2021fourier}, flow map learning
(FML)~\cite{churchill2023flow}, and neural operators defined on POD bases~\cite{Bhattacharya_2021}. These methods differ in the types of
operators they approximate well.
We employ FML because of our focus on predicting the level-set time
evolution, along with the use of with neural operators defined on POD bases
to facilitate learning from a limited number of training simulations. Recent work has extended
FML from its original focus on autonomous dynamical
systems~\cite{qin2020data} to include control
inputs~\cite{qin2021data}, parametric dependence~\cite{qin2021deep},
modal-space representations such as POD or Fourier
series~\cite{wu2020data, chen2022deep, churchill2025principal}, and
flow-map-constrained optimization~\cite{hart2023solving}. 

In FML, it is frequently advantageous to learn the ROM
on a coarse time grid since this can reduce the effect of error accumulation. Accordingly, we assume that the WRF-SFIRE simulation data is subsampled at $n_\tau=9$ points that correspond to the times $t=0,1,\dots,8$ hours. Letting $\mathring{\u}(\z)$ denote the spatial discretization of the initial condition parameterization~\eqref{eqn:level_set_init}, we consider the discrete-time dynamical system
\begin{equation}
    \begin{rcases}
        \u_0 &= \mathring{\u} (\z) \\
        \u_{k} &= \u_{k-1} + \f(\u_{k-1}), \qquad k = 1, \dots, n_\tau-1 \qquad
    \end{rcases}
    \label{equ:wrf-flowmap-full}
\end{equation}
where $\f:\R^{n_u} \to \R^{n_u}$ corresponds to the integrating the
WRF-SFIRE dynamics over a one-hour time interval.  Note that $\f$ also depends on the atmospheric winds, although it is not
explicit in the notation. Learning the coupling between the level-set function and the wind field from limited training simulations is challenging due to the complex nonlinear interactions between fire and the atmosphere. Instead, we will fix the atmospheric wind initial condition to a nominal estimate for all training data generation and only extract the level-set function as output from the WRF-SFIRE simulations. This corresponds, in our previous notation, to fixing $\xivec=\tilde{\xivec}$, where $\xivec$ represents the atmospheric wind initialization. Due to the fire spread model's dependence on the transient wind field,~\eqref{equ:wrf-flowmap-full} cannot represent the WRF-SFIRE dynamics exactly. However, with a fixed wind field initialization,~\eqref{equ:wrf-flowmap-full} is a reasonable approximation to facilitate learning a ROM.

\begin{figure}
    \centering
    \includegraphics[width=0.99\textwidth]{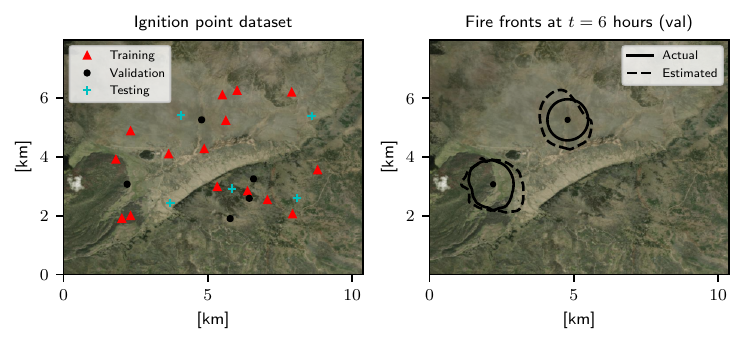}
    \caption{Left: Satellite image of the Valles Caldera with 25 ignition locations, divided into training, validation, and test sets. Right: Actual and estimated fire fronts from the validation set.}
    \label{fig:valles-terrain-locs}
\end{figure}

 An initial wind field $\v(\x, 0) = (0.1, 0.1, 0)^\top$ is fixed and the 20 simulations are run with different ignition locations to generate a training set (with 15 simulations) and validation set (with 5 simulations). The left panel of~\Cref{fig:valles-terrain-locs} shows the terrain and ignition locations.

Our objective in FML is to approximate $\f$ in~\eqref{equ:wrf-flowmap-full} using our training and validation set simulations. Since the dimension of a state time snapshot is $\mathcal{O}(10^4)$, it is impractical to learn $\f:\R^{n_u} \to \R^{n_u}$ from only $20$ simulations. Rather, we project the dynamics into a POD basis and approximate the projected flow map. At least $r=100$ POD basis vectors are required to achieve a $1\%$ relative $\ell^2$ projection error on the validation trajectories.

Let $\V \in \R^{n_u \times r}$ denote the matrix containing the POD basis vectors. Multiplying~\cref{equ:wrf-flowmap-full} by $\V^\top$ on the left motivates the FML-based reduced dynamical system
\begin{align}\label{eqn:reduced_flow_map}
    \begin{rcases}
        \widehat{\u}_0 &=  \V^\top \mathring{\u} (\vec{z}) \\
        \widehat{\u}_{k} &= \widehat{\u}_{k-1} + \widehat{\f}(\widehat{\u}_{k-1}) \qquad k=1, \dots, n_\tau-1 \qquad
    \end{rcases} \ ,
\end{align}
where $\widehat{\f}(\V^\top \u_k)$ is a function we seek to learn that approximates the reduced WRF-SFIRE dynamics $\V^\top \f(\V \V^\top \u_k)$. Appendix~\ref{appendix:wildfire} provides details regarding the network architecture of $\hat{\f}$, loss function, and training. The right panel of~\Cref{fig:valles-terrain-locs} displays the fire front contour at $t=6$ hours for a validation set simulation and the corresponding ROM prediction.

\subsection{ROM-constrained optimization}

Once we have trained the ROM, the \texttt{FOMCO} problem~\cref{equ:wrf-objective} is approximated using the \texttt{ROMCO} problem
\begin{equation}
    \begin{array}{cl}
        \displaystyle \min_{\z \in \R^2}  & \displaystyle \sum_{j=1}^{n_\text{obs}} \|  \V \widehat{\u}_{1+j}(\z) - \y_j \|^2
    \end{array} \ \
    \label{equ:wrf-poco-opt}
\end{equation}
where $\widehat{\u}_{1+j}(\z)$, $j=1,2,\dots,n_\text{obs}$, satisfies~\cref{eqn:reduced_flow_map}. To assess how well~\eqref{equ:wrf-poco-opt} recovers the ignition location, we generate a test set that consists of $5$ WRF-SFIRE simulations with ignition locations not used in the training and validation sets (\Cref{fig:valles-terrain-locs} depicts the ignition locations). Furthermore, the test set simulations use the initial wind field $\v(\x, 0) = (4.33, 2.5, 0)^\top$, which has a significantly different magnitude and direction than the training and validation-set initial-wind fields. In the article's notation, this test set wind field initialization corresponds to $\xivec^\dagger$. 

For each testing dataset simulation, we consider the level-set from WRF-SFIRE as noise-free observations $\{ \y_j \}_{j=1}^{n_\text{obs}}$ that enter the optimization objective. A trust-region method is used to solve~\eqref{equ:wrf-poco-opt} for each test set simulation. Gradients are computed by using automatic differentiation
to compute the Jacobian of $\widehat{\f}$, and the Hessian
is approximated using BFGS~\cite{nocedal_wright_book}.
The data
misfit---the difference between the test set data and the FOM solution
using the estimated ignition location---is shown in the second column
of~\Cref{tbl:valles-objective-vals}. The rows correspond to each of the $5$
test cases. We observe misfit errors ranging from
$0.1\%$ to $19\%$, with an average of $5.4\%$. The ignition location
estimation error---the difference between the true and estimated
ignition location---is shown in the second column
of~\Cref{tbl:valles-z-vals}. The ignition location errors range from
$281$ to $2145$ meters over the $5$ test cases, with an average of $967$ meters.
This ignition location error is what we seek to reduce via HDSA-MD, which is described below and reported in the right columns of the tables.

\subsection{Model discrepancy analysis}

The large errors we observe in~\Cref{tbl:valles-objective-vals} and~\Cref{tbl:valles-z-vals} highlight the deficiency of using a ROM trained with inaccurate atmospheric data. Nonetheless, each evaluation of WRF-SFIRE consists of taking over $10^4$ time steps, which requires considerable wall-clock-time, making it intractable to solve the \texttt{FOMCO} problem in real-time. Our goal is to get the best ignition estimate possible using only one WRF-SFIRE evaluation online.

For each of the 5 test cases, we generate a single FOM simulation
using the test-set initial-wind field $\v(\x, 0) = (4.33, 2.5, 0)^\top$
and ignition location given by the \texttt{ROMCO} solution $\ztilde$
for the given test case. This single FOM evaluation is used to calibrate the model discrepancy. Let $\overline{\t}$ denote the posterior mean of the discrepancy calibration and $\overline{\z}=\tilde{\z}-\H^{-1} \B \overline{\t}$ denote the posterior mean ignition location. For each of the 5 test cases, we run WRF-SFIRE with the ignition location $\overline{\z}$ and test set initial wind field to measure the reduction in data misfit. The third columns of~\Cref{tbl:valles-objective-vals} and~\Cref{tbl:valles-z-vals} display the posterior mean ignition location data misfit and ignition location estimation error, respectively, for each of the test cases. The fourth column of~\Cref{tbl:valles-objective-vals} and~\Cref{tbl:valles-z-vals} shows the relative improvement from the \texttt{ROMCO} solution to the posterior mean ignition location. We observe $36\% - 74\%$ reduction in the data misfit and $18\% - 59\%$ reduction in the ignition location estimation error. The lowest percentage reduction corresponds to the test case where the error in $\tilde{\z}$ is already small. The consistent improvement attained using a single FOM evaluation in each test case demonstrates the benefit of HDSA-MD.

\begin{table}[t]
    \centering
    \begin{tabular}{c||c|c|c}
        Test Case & Relative misfit at $\tilde{\z}$ & Relative misfit at $\overline{\z}$ &  Reduction \\ \hline
        1 & $1.25 \times 10^{-3}$ & $7.98 \times 10^{-4}$ & 36\% \\
        2 & $1.90 \times 10^{-1}$ & $8.23 \times 10^{-2}$ & 57\% \\
        3 & $1.74 \times 10^{-3}$ & $4.51 \times 10^{-4}$ & 74\% \\
        4 & $3.67 \times 10^{-2}$ & $1.22 \times 10^{-2}$ & 67\% \\
        5 & $4.13 \times 10^{-2}$ & $1.71 \times 10^{-2}$ & 59\% \\ \hline
        \textbf{Average} & $5.42 \times 10^{-2}$ & $2.26 \times 10^{-2}$ & 58\%
    \end{tabular}
    \caption{FOM misfit (normalized by $\| \vec{Y} \|^2_{F}$) from running WRF-SFIRE at the ROMCO solution $\tilde{\z}$ and the posterior mean ignition location $\overline{\z}=\tilde{\z}-\H^{-1} \B \overline{\t}$.}
    \label{tbl:valles-objective-vals}
\end{table}

\begin{table}[t]
    \centering
    \begin{tabular}{c||c|c|c}
        Test Case & $\| \tilde{\z}- \z^\dagger \|$ & $\| \overline{\z} - \z^\dagger \|$ & Reduction \\ \hline
        1 & $2.81 \times 10^{2}$ & $2.31 \times 10^{2}$ & 18\% \\
        2 & $2.14 \times 10^{3}$ & $1.41 \times 10^{3}$ & 34\% \\
        3 & $3.26 \times 10^{2}$ & $1.33 \times 10^{2}$ & 59\% \\
        4 & $9.99 \times 10^{2}$ & $5.45 \times 10^{2}$ & 45\% \\
        5 & $1.08 \times 10^{3}$ & $6.60 \times 10^{2}$ & 39\% \\  \hline
        \textbf{Average} & $9.67 \times 10^{2}$ & $5.96 \times 10^{2}$ & 38\%
    \end{tabular}
    \caption{Error (in meters) of the ROMCO solution $\tilde{\z}$ and the posterior mean ignition location $\overline{\z}$; $\z^\dagger$ denotes the ``ground truth'' ignition location used in the test set simulation.}
    \label{tbl:valles-z-vals}
\end{table}

\begin{figure}[h]
    \centering
    \includegraphics[width=0.49\textwidth]{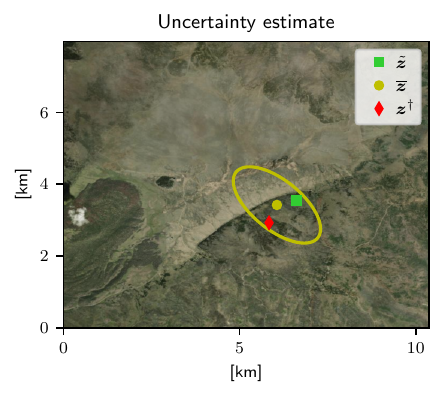}
    \caption{Posterior mean and uncertainty in the optimization solution for a single test case. The dashed line defines the 95\% confidence region.}
    \label{fig:valles-optima-and-updates}
\end{figure}

Lastly, we select one test case to demonstrate the utility of the
model discrepancy analysis to provide uncertainty
estimates. \Cref{fig:valles-optima-and-updates} displays the
posterior mean ignition location $\overline{\z}$ and a $95 \%$ confidence ellipse around it (the posterior follows a Gaussian distribution), along with the
\texttt{ROMCO} solution $\tilde{\z}$ and the ``ground truth'' ignition
location $\z^\dagger$. We observe that the posterior
 mean ignition location lies approximately half-way between the \texttt{ROMCO} solution
 and the true ignition location, reflecting the approximately $50\%$
 reduction in estimation error. Furthermore, the true ignition
location $\z^\dagger$ lies within the $95 \%$ confidence ellipse.
Hence, the posterior distribution is well-calibrated despite the fact that it
only utilized a single FOM evaluation.

\section{Conclusion} \label{sec:conclusion}

Enabling optimization within wall-clock-time constraints is essential for developing digital twins that adapt to incoming data and inform decision-making. This article introduces a novel combination of non-intrusive ROMs and HDSA-MD to address the challenge of rapidly computing a reliable optimization solution. Our numerical results demonstrate ROM-constrained optimization for two different systems of nonlinear partial differential equations. In both cases, the ROMs are trained using inaccurate atmospheric data, which is subsequently corrected using HDSA-MD with a single FOM evaluation involving updated atmospheric data. This is representative of a range of digital-twin applications where model parameters vary with time, and hence the ROM is necessarily trained with inaccurate parameters. 

In scenarios where it is feasible, retraining the ROM using the additional FOM evaluations and resolving the updated optimization problem may be superior. Our approach is closely related since HDSA-MD consists of (i) calibrating a model discrepancy representation to augment the ROM, and (ii) performing a post-optimality update to approximate the discrepancy augmented optimization solution. The advantage of HDSA-MD is that it avoids parameter tuning and additional computation required for ROM retraining, thus making the process automatic and real-time performant. Furthermore, HDSA-MD is effective with as little as one additional FOM evaluation. Lastly, HDSA-MD provides efficient estimates of uncertainty due to the model discrepancy.

\section*{Acknowledgements}

This work was supported by the U.S. Department of Energy,
Office of Science, Office of Advanced Scientific Computing Research Field Work
Proposal Number 23-02526, and the Laboratory Directed
Research and Development program at Sandia National Laboratories. Sandia National Laboratories, a
multimission laboratory managed and operated by National Technology and Engineering Solutions of
Sandia LLC, a wholly owned subsidiary of Honeywell International Inc. for the U.S. Department of Energy’s
National Nuclear Security Administration contract DE-NA0003525.
This written work is authored by an employee of NTESS. The employee, not NTESS, owns the right, title and interest in and to the written work and is responsible for its contents. Any subjective views or opinions that might be expressed in the written work do not necessarily represent the views of the U.S. Government. The publisher acknowledges that the U.S. Government retains a non-exclusive, paid-up, irrevocable, world-wide license to publish or reproduce the published form of this written work or allow others to do so, for U.S. Government purposes. The DOE will provide public access to results of federally sponsored research in accordance with the DOE Public Access Plan. SAND2025-10864O.

%% file: appendices.tex
\section{Operator Inference ROM for Atmospheric Contaminant Modeling} \label{appendix:atmosphere_transport}

\paragraph{Training data generation} We construct $n_c = 5$ smooth random control vectors $\q^{(\ell)}(t)$, $\ell=1,\ldots,n_c$, by drawing values for $\q^{(\ell)}(\bar{t})$ at $\bar{t} = 0$, $T/3$, $2T/3$, and $T$ from the uniform distribution over $[0, 5]$ and interpolating these values with a monotone piecewise cubic spline~\cite{fritsch1980pchipinterp}.
The full-order system~\cref{eq:adr-ode} is then solved with $\q(t;\z) = \q^{(\ell)}(t)$, resulting in snapshot matrices
$\U_1^{(\ell)},\U_2^{(\ell)}\in\R^{n_x \times n_t}$ for each $\ell=1,\ldots,n_c$. Here, the $j$-th column of $\U_k^{(\ell)}$ is the $\u_k(t_j)$ (the $k$-th variable at $j$-th time) with forcing $\q^{(\ell)}$.

\paragraph{POD basis generation} The basis matrices $\V_1,\V_2$ are obtained through weighted proper orthogonal decomposition (POD), which minimizes the reconstruction error of the snapshot data in the $\M$-norm \cite{algazi1969karhunenloeve,berkooz1993pod,sirovich1987pod}.
Specifically, for $k = 1, 2$, we form the concatenated snapshot matrix
$\U_k = (\U_k^{(1)}~\U_k^{(2)}~\cdots~\U_k^{(n_c)})\in\R^{n_x \times n_c n_t}$,
compute the symmetric eigendecomposition $\U_k\trp\M\U_k = \vec{\Psi}_k\vec{\Sigma}_k^2\vec{\Psi}_k\trp$, and define $\V_k$ to be the first $r_k$ columns of $\U_k\vec{\Psi}_k\vec{\Sigma}_k^{-1}$.
The reduced dimension $r_k$ is chosen based on the POD singular values---the diagonal values of $\vec{\Sigma}_k$, written $\sigma_1>\sigma_2>\cdots>\sigma_{n_c n_t} > 0$---so that the so-called POD residual energy is less than $10^{-5}$, i.e.,
\begin{align}
    \label{eq:adr-podres}
    \left.\sum_{j=r_k+1}^{n_c n_t}\sigma_j^2\, \middle/ ~\sum_{j=1}^{n_c n_t}\sigma_j^2\right. < 10^{-5}.
\end{align}
In our experiment, the smallest $r_1,r_2$ that satisfy this criterion are $r_1 = 19$ and $r_2 = 33$.

\paragraph{Operator inference regression problem}  To generate the ROM, we formulate a linear regression problem to estimate the operators~\cref{eq:adr-galerkin-matrices}. To this end, we project the snapshot matrices onto the bases $\V_k$, $k=1,2$, i.e., $(\widehat{\u}_k^1~\cdots~\widehat{\u}_k^{K}) = \V_k\trp\M\U_k$, with the number of columns $K = n_c n_t$ corresponding the product of the number of FOM evaluations $n_c$ and the number of time snapshots $n_t$. The operator inference regression problem requires access to time derivatives, which we denote as $\dot{\widehat{\u}}_k^j$. We approximate the time derivatives using sixth-order finite differences of the training states. Then we estimate the operators by independently solving the two linear regression problems
\begin{subequations}
\label{eq:adr-opinf}
\begin{align}
    \label{eq:adr-opinf-u1}
    \min_{\widehat{\A}_1,\widehat{\mathbf{R}}_1}\sum_{j=1}^{K}\left\|
        \widehat{\A}_1\widehat{\u}_1^j
        + \widehat{\mathbf{R}}_1 \left( \widehat{\u}_1^j\otimes\widehat{\u}_2^j \right)
        - \dot{\widehat{\u}}_1^j
    \right\|_{2}^{2}
    &+ \lambda_1\|\widehat{\A}_1\|_F^2
    + \lambda_2\|\widehat{\mathbf{R}}_1\|_F^2,
    \\
    \label{eq:adr-opinf-u2}
    \min_{\widehat{\A}_2,\widehat{\mathbf{R}}_2,\widehat{\vec{\Phi}}}\sum_{j=1}^{K}\left\|
        \widehat{\A}_2\widehat{\u}_2^j
        + \widehat{\mathbf{R}}_2 \left( \widehat{\u}_1^j\otimes\widehat{\u}_2^j \right)
        + \widehat{\vec{\Phi}}\left( \q^j\odot\q^j \right)
        - \dot{\widehat{\u}}_1^j
    \right\|_{2}^{2}
    &+ \lambda_1\left(\|\widehat{\A}_2\|_F^2 + \|\widehat{\vec{\Phi}}\|_F^2\right)
    + \lambda_2\|\widehat{\mathbf{R}}_2\|_F^2.
\end{align}
\end{subequations}
To select appropriate regularization constants $\lambda_1,\lambda_2 \ge 0$, we take the approach of \cite{mcquarrie2021regopinf,mcquarrie2023popinf,qian2022pdeopinf}: for each pair in a set of potential values $(\lambda_1,\lambda_2)$, we solve the operator inference regression problem~\cref{eq:adr-opinf} and evaluate the resulting ROM~\cref{eq:adr-rom} with the same initial conditions and controls as were used to generate the full-order data, then select the $(\lambda_1,\lambda_2)$ that minimizes the average reconstruction error.
In this example, the process results in $\lambda_1 = 0.1$ and $\lambda_2 = 10$.

\paragraph{Remark on efficient objective function evaluation} The integrand in the objective \cref{eq:adr-objective-rom} can be evaluated efficiently by pre-computing norm weights for the reduced state.
Let $\tilde{\M}\in\R^{n_x\times n_x}$ be the matrix with entries $\tilde{M}_{ij} = \psi_i M_{ij} \psi_j$, where the $M_{ij}$ are the entries of the mass matrix $\M$ and the $\psi_i$'s are the entries of $\vec{\psi}$.
Then
\begin{align}
    \|\V_1\widehat{\u}_1(t) \odot \vec{\psi} \|_{\M}^2
    = (\V_1\widehat{\u}_1(t) \odot \vec{\psi} )\trp\M(\V_1\widehat{\u}_1(t) \odot \vec{\psi})
    = \widehat{\u}_1(t)\trp\V_1\trp\tilde{\M}\V_1\widehat{\u}_1(t).
\end{align}
The matrix $\V_1\trp\tilde{\M}\V_1\in\R^{r_1\times r_1}$ can be computed once and stored for repeated use so that the cost of evaluating the objective scales with $r_1$ instead of $n_x$.

\section{Neural operator ROM for Wildfire Modeling} \label{appendix:wildfire}

We seek a neural network approximation $\widehat{\f}(\V^\top \u \, ; \vec{\gamma}) \approx \V^\top \f (\V \V^\top \u)$, where $\vec{\gamma}$ are the network parameters (weights and biases). Let $\widehat{\u}$ denote the approximate reduced state coordinates generated by repeated composition of $\widehat{\f}$, as in~\cref{eqn:reduced_flow_map}, where $\widehat{\u}$ and $\widehat{\f}$ depend parametrically on $\vec{\gamma}$.

Let $M=16$ denote the number of training dataset simulations and $\u^{(m)}_k$ denote the $m^{th}$ training dataset state solution at time $t_k$, $m=1,2,\dots,M$, $k=0,1,\dots,n_\tau$. To train $\widehat{\f}$, we seek to determine weights and biases $\vec{\gamma}$ that minimize the $P$-step recurrent loss
\begin{equation}
    L(\vec{\gamma}) = \sum_{m=1}^M \, \sum_{k=0}^{n_\tau} \, \sum_{p=1}^P \, \left\|  \widehat{\u}^{(m)}_{k+p}(\vec{\gamma}) - \V^\top \u^{(m)}_{k+p} \right\|_{\ell^2}^2,
    \label{equ:wrf-flowmap-loss}
\end{equation}
where the innermost sum is assumed to be zero if $k+p > n_\tau$ and $\V$ is the POD basis matrix. Setting $1 < P < n_\tau$ seeks to balance the long-time stability of $\widehat{\u}$ against the additional training costs of incorporating further compositions of $\widehat{\f}$~\cite{hart2023solving, patel2021physics}.

The flow-map network has 4 hidden layers with 200 neurons per layer and a GELU activation function between layers; the output layer has an identity activation. Training is done with Tensorflow using the Adam optimizer. We train for 1000 epochs with an initial learning rate of $4 \times 10^{-3}$ that decays smoothly to $10^{-3}$. At the end of training, the flow-map has a relative $\ell^2$ composition error of $2.3\%$ over the validation set.

%% file: main.bbl
\begin{thebibliography}{10}

\bibitem{algazi1969karhunenloeve}
V.~Algazi and D.~Sakrison.
\newblock On the optimality of the {K}arhunen--{L}o{\`e}ve expansion.
\newblock {\em IEEE Transactions on Information Theory}, 15(2):319--321, 1969.

\bibitem{antil2011DDandBT}
H.~Antil, M.~Heinkenschloss, and R.~Hoppe.
\newblock Domain decomposition and balanced truncation model reduction for
  shape optimization of the {S}tokes system.
\newblock {\em Optimization Methods and Software}, 26(4-5), 2011.

\bibitem{antil2012ROM}
H.~Antil, M.~Heinkenschloss, R.~Hoppe, C.~Linsenmann, and A.~Wixforth.
\newblock Reduced order modeling based shape optimization of surface acoustic
  wave driven microfluidic biochips.
\newblock {\em Mathematics and Computers in Simulation}, 82(10), 2012.

\bibitem{berkooz1993pod}
G.~Berkooz, P.~Holmes, and J.~L. Lumley.
\newblock The proper orthogonal decomposition in the analysis of turbulent
  flows.
\newblock {\em Annual Review of Fluid Mechanics}, 25:539--575, 1993.

\bibitem{Bhattacharya_2021}
K.~Bhattacharya, B.~Hosseini, N.~B. Kovachki, and A.~M. Stuart.
\newblock Model reduction and neural networks for {PDEs}.
\newblock {\em SMAI Journal of Computational Mathematics}, 7:121--157, 2021.

\bibitem{ghattas_infinite_dim_bayes_1}
T.~Bui-{T}hanh, O.~Ghattas, J.~Martin, and G.~Stadler.
\newblock A computational framework for infinite-dimensional {B}ayesian inverse
  problems. {P}art {I}: The linearized case, with applications to global
  seismic inversion.
\newblock {\em SIAM Journal on Scientific Computing}, 35(6):A2494--A2523, 2013.

\bibitem{buithanh2006goob}
T.~Bui-Thanh, K.~Willcox, O.~Ghattas, and B.~{van Bloemen Waanders}.
\newblock Goal-oriented, model-constrained optimization for reduction of
  large-scale systems.
\newblock {\em Journal of Computational Physics}, 224(2):880--896, 2007.

\bibitem{chen2022deep}
Z.~Chen, V.~Churchill, K.~Wu, and D.~Xiu.
\newblock Deep neural network modeling of unknown partial differential
  equations in nodal space.
\newblock {\em Journal of Computational Physics}, 449:110782, 2022.

\bibitem{churchill2025principal}
V.~Churchill.
\newblock Principal component flow map learning of {PDE}s from incomplete,
  limited, and noisy data.
\newblock {\em Journal of Computational Physics}, 524:113730, 2025.

\bibitem{churchill2023flow}
V.~Churchill and D.~Xiu.
\newblock Flow map learning for unknown dynamical systems: overview,
  implementation, and benchmarks.
\newblock {\em Journal of Machine Learning for Modeling and Computing},
  4(2):173--201, 2023.

\bibitem{fritsch1980pchipinterp}
F.~N. Fritsch and R.~E. Carlson.
\newblock Monotone piecewise cubic interpolation.
\newblock {\em SIAM Journal on Numerical Analysis}, 17(2):238--246, 1980.

\bibitem{ghattas2021physicsbased}
O.~Ghattas and K.~Willcox.
\newblock Learning physics-based models from data: perspectives from inverse
  problems and model reduction.
\newblock {\em Acta Numerica}, 30:445--554, 2021.

\bibitem{hassdonk2023rbmlromopt}
B.~Haasdonk, H.~Kleikamp, M.~Ohlberger, F.~Schindler, and T.~Wenzel.
\newblock A new certified hierarchical and adaptive {RB}-{ML}-{ROM} surrogate
  model for parametrized {PDE}s.
\newblock {\em SIAM Journal on Scientific Computing}, 45(3):A1039--A1065, 2023.

\bibitem{hart2023solving}
J.~Hart, M.~Gulian, I.~Manickam, and L.~P. Swiler.
\newblock Solving high-dimensional inverse problems with auxiliary uncertainty
  via operator learning with limited data.
\newblock {\em Journal of Machine Learning for Modeling and Computing},
  4(2):105--133, 2023.

\bibitem{hart2023hdsaoptimal}
J.~Hart and B.~{van Bloemen Waanders}.
\newblock Hyper-differential sensitivity analysis with respect to model
  discrepancy: {O}ptimal solution updating.
\newblock {\em Computer Methods in Applied Mechanics and Engineering}, 412,
  2023.

\bibitem{hart2024hdsasampling}
J.~Hart and B.~{van Bloemen Waanders}.
\newblock Hyper-differential sensitivity analysis with respect to model
  discrepancy: {P}osterior optimal solution sampling.
\newblock {\em {AIMS} Foundations of Data Science}, 2024.

\bibitem{hart_vbw_2025}
J.~Hart, B.~{van Bloemen Waanders}, J.~Li, T.~A. Ouermi, and C.~R. Johnson.
\newblock Hyper-differential sensitivity analysis with respect to model
  discrepancy: Prior distributions.
\newblock {\em Under review. arxiv/2504.19812}, 2025.

\bibitem{heinkenschloss2018rom}
M.~Heinkenschloss and D.~Jando.
\newblock Reduced order modeling for time-dependent optimization problems with
  initial value controls.
\newblock {\em SIAM Journal on Scientific Computing}, 40(1):A22--A51, 2018.

\bibitem{heinkenschloss2018CVaR}
M.~Heinkenschloss, B.~Kramer, T.~Takhtaganov, and K.~Willcox.
\newblock Conditional-value-at-risk estimation via reduced-order models.
\newblock {\em SIAM Journal on Uncertainty Quantification}, 6(4):1395--1423,
  2018.

\bibitem{kovachki2023neural}
N.~Kovachki, Z.~Li, B.~Liu, K.~Azizzadenesheli, K.~Bhattacharya, A.~Stuart, and
  A.~Anandkumar.
\newblock Neural operator: Learning maps between function spaces with
  applications to {PDEs}.
\newblock {\em Journal of Machine Learning Research}, 24(89):1--97, 2023.

\bibitem{kramer2024opinfsurvey}
B.~Kramer, B.~Peherstorfer, and K.~E. Willcox.
\newblock Learning nonlinear reduced models from data with operator inference.
\newblock {\em Annual Review of Fluid Mechanics}, 56(1):521--548, 2024.

\bibitem{li2021fourier}
Z.~Li, N.~B. Kovachki, K.~Azizzadenesheli, B.~Liu, K.~Bhattacharya, A.~Stuart,
  and A.~Anandkumar.
\newblock Fourier neural operator for parametric partial differential
  equations.
\newblock In {\em International Conference on Learning Representations}, 2021.

\bibitem{lu2021learning}
L.~Lu, P.~Jin, G.~Pang, Z.~Zhang, and G.~E. Karniadakis.
\newblock Learning nonlinear operators via {DeepONet} based on the universal
  approximation theorem of operators.
\newblock {\em Nature Machine Intelligence}, 3:218--229, 2021.

\bibitem{Luo_2025}
D.~Luo, T.~O'Leary-Roseberry, P.~Chen, and O.~Ghattas.
\newblock Efficient {PDE}-constrained optimization under high-dimensional
  uncertainty using derivative-informed neural operators.
\newblock {\em SIAM Journal on Scientific Computing}, 47(4):C899--C931, 2025.

\bibitem{mandel2011coupled}
J.~Mandel, J.~D. Beezley, and A.~K. Kochanski.
\newblock Coupled atmosphere--wildland fire modeling with {WRF} 3.3 and {SFIRE}
  2011.
\newblock {\em Geoscientific Model Development}, 4(3):591--610, 2011.

\bibitem{mcquarrie2021regopinf}
S.~A. McQuarrie, C.~Huang, and K.~E. Willcox.
\newblock Data-driven reduced-order models via regularised operator inference
  for a single-injector combustion process.
\newblock {\em Journal of the Royal Society of New Zealand}, 51(2):194--211,
  2021.

\bibitem{mcquarrie2023popinf}
S.~A. McQuarrie, P.~Khodabakhshi, and K.~E. Willcox.
\newblock Non-intrusive reduced-order models for parametric partial
  differential equations via data-driven operator inference.
\newblock {\em SIAM Journal on Scientific Computing}, 45(4):A1917--A1946, 2023.

\bibitem{NAP26894}
{National Academies of Sciences, Engineering, and Medicine}.
\newblock {\em Foundational Research Gaps and Future Directions for Digital
  Twins}.
\newblock The National Academies Press, Washington, DC, 2024.

\bibitem{nocedal_wright_book}
J.~Nocedal and S.~J. Wright.
\newblock {\em Numerical {O}ptimization}.
\newblock Springer, New York, 2nd edition, 2006.

\bibitem{dipnet_2022}
T.~O'Leary-Roseberry, U.~Villa, P.~Chen, and O.~Ghattas.
\newblock Derivative-informed projected neural networks for high-dimensional
  parametric maps governed by {PDE}s.
\newblock {\em Computer Methods in Applied Mechanics and Engineering},
  388:114199, 2022.

\bibitem{patel2021physics}
R.~G. Patel, N.~A. Trask, M.~A. Wood, and E.~C. Cyr.
\newblock A physics-informed operator regression framework for extracting
  data-driven continuum models.
\newblock {\em Computer Methods in Applied Mechanics and Engineering},
  373:113500, 2021.

\bibitem{peherstorfer2016opinf}
B.~Peherstorfer and K.~Willcox.
\newblock Data-driven operator inference for nonintrusive projection-based
  model reduction.
\newblock {\em Computer Methods in Applied Mechanics and Engineering},
  306:196--215, 2016.

\bibitem{qian2022pdeopinf}
E.~Qian, I.-G. Farca\c{s}, and K.~Willcox.
\newblock Reduced operator inference for nonlinear partial differential
  equations.
\newblock {\em SIAM Journal on Scientific Computing}, 44(4):A1934--A1959, 2022.

\bibitem{qian2017rbopt}
E.~Qian, M.~Grepl, K.~Veroy, and K.~Willcox.
\newblock A certified trust region reduced basis approach to {PDE}-constrained
  optimization.
\newblock {\em SIAM Journal on Scientific Computing}, 39(5):S434--S460, 2017.

\bibitem{qin2021data}
T.~Qin, Z.~Chen, J.~D. Jakeman, and D.~Xiu.
\newblock Data-driven learning of nonautonomous systems.
\newblock {\em SIAM Journal on Scientific Computing}, 43(3):A1607--A1624, 2021.

\bibitem{qin2021deep}
T.~Qin, Z.~Chen, J.~D. Jakeman, and D.~Xiu.
\newblock Deep learning of parameterized equations with applications to
  uncertainty quantification.
\newblock {\em International Journal for Uncertainty Quantification},
  11(2):63--82, 2021.

\bibitem{qin2020data}
T.~Qin, K.~Wu, and D.~Xiu.
\newblock Data driven governing equations approximation using deep neural
  networks.
\newblock {\em Journal of Computational Physics}, 395:620--635, 2020.

\bibitem{shaddy2024generative}
B.~Shaddy, D.~Ray, A.~Farguell, V.~Calaza, J.~Mandel, J.~Haley, K.~Hilburn,
  D.~V. Mallia, A.~Kochanski, and A.~Oberai.
\newblock Generative algorithms for fusion of physics-based wildfire spread
  models with satellite data for initializing wildfire forecasts.
\newblock {\em Artificial Intelligence for the Earth Systems}, 3(3):e230087,
  2024.

\bibitem{sirovich1987pod}
L.~Sirovich.
\newblock Turbulence and the dynamics of coherent structures. {I}. {C}oherent
  structures.
\newblock {\em Quarterly of Applied Mathematics}, 45(3):561--571, 1987.

\bibitem{skamarock2019wrf}
W.~Skamarock, J.~Klemp, J.~Dudhia, D.~O. Gill, Z.~Liu, J.~Berner, W.~Wang,
  J.~G. Powers, M.~G. Duda, D.~Barker, and X.-Y. Huang.
\newblock A description of the advanced research {WRF} model version 4.1.
\newblock Technical report, National Center for Atmospheric Research, 2019.

\bibitem{stuart_inv_prob}
A.~M. Stuart.
\newblock Inverse problems: {A} {B}ayesian perspective.
\newblock {\em Acta Numerica}, 19:451--559, 2010.

\bibitem{wu2020data}
K.~Wu and D.~Xiu.
\newblock Data-driven deep learning of partial differential equations in modal
  space.
\newblock {\em Journal of Computational Physics}, 408:109307, 2020.

\bibitem{yue2013optlinearrom}
Y.~Yue and K.~Meerbergen.
\newblock Accelerating optimization of parametric linear systems by model order
  reduction.
\newblock {\em SIAM Journal on Optimization}, 23(2):1344--1370, 2013.

\bibitem{zahr2015pderom}
M.~J. Zahr and C.~Farhat.
\newblock Progressive construction of a parametric reduced-order model for
  {PDE}-constrained optimization.
\newblock {\em International Journal for Numerical Methods in Engineering},
  102(5):1111--1135, 2015.

\bibitem{zou2022CVaR}
Z.~Zou, D.~P. Kouri, and W.~Aquino.
\newblock A locally adapted reduced-basis method for solving risk-averse
  {PDE}-constrained optimization problems.
\newblock {\em SIAM Journal on Uncertainty Quantification}, 10(4), 2022.

\end{thebibliography}
